\documentclass{amsart}
\usepackage[utf8]{inputenc}

\usepackage{textcmds}
\usepackage{mathtools}
\usepackage{amssymb}
\usepackage{amsthm}
\usepackage{amsmath}
\usepackage{bm}
\usepackage{bbm}
\usepackage{enumerate}
\usepackage{cite}
\usepackage{xfrac}
\usepackage{xcolor}
\usepackage{tikz-cd}
\usepackage{lipsum}
\usepackage{adjustbox}
\usepackage{url}
\usepackage{hyperref}
\usepackage{multicol}

\usepackage[left=3cm,right=3cm,top=3cm,bottom=4.5cm]{geometry}

\theoremstyle{plain}
\newtheorem{theorem}{Theorem}[section]

\newtheorem{lemma}[theorem]{Lemma}
\newtheorem*{lemma*}{Lemma}

\theoremstyle{definition}
\newtheorem{definition}[theorem]{Definition}
\newtheorem*{definition*}{Definition}

\theoremstyle{plain}

\theoremstyle{plain}

\theoremstyle{plain}

\theoremstyle{plain}
\newtheorem{fact}[theorem]{Fact}

\theoremstyle{remark}
\newtheorem{remark}[theorem]{Remark}

\theoremstyle{remark}
\newtheorem{notation}[theorem]{Notation}

\theoremstyle{remark}

\theoremstyle{remark}
\newtheorem{convention}[theorem]{Convention}

\theoremstyle{plain}
\newtheorem{proposition}[theorem]{Proposition}

\theoremstyle{plain}

\newtheorem{theoremintro}{Theorem}

\theoremstyle{definition}

\def\dotminussym#1{%
  \setbox0=\hbox{$-$}%
  \kern.5\wd0%
  \hbox to 0pt{\hss\hbox{$-$}\hss}%
  \raise.6\ht0\hbox to 0pt{\hss$.$\hss}%
  \kern.5\wd0%
}

\def\Ind#1#2{#1\setbox0=\hbox{$#1x$}\kern\wd0\hbox to 0pt{\hss$#1\mid$\hss}
\lower.9\ht0\hbox to 0pt{\hss$#1\smile$\hss}\kern\wd0}

\def\ind{\mathop{\mathpalette\Ind{}}}

\def\notind#1#2{#1\setbox0=\hbox{$#1x$}\kern\wd0
\hbox to 0pt{\mathchardef\nn=12854\hss$#1\nn$\kern1.4\wd0\hss}
\hbox to 0pt{\hss$#1\mid$\hss}\lower.9\ht0 \hbox to 0pt{\hss$#1\smile$\hss}\kern\wd0}

\title{Higher amalgamation in $\mathrm{ACFA}^{+}$}

\author{Stefan Marian Ludwig}
\address{Albert-Ludwigs-Universität Freiburg,
Mathematisches Institut,
Abteilung für Mathematische Logik,
Ernst-Zermelo-Straße 1,
79104 Freiburg i. B., Germany}
\email{stefan.ludwig@mathematik.uni-freiburg.de}	
\thanks{SML has received funding from the European Union's Horizon 2020 research and innovation programme under the Marie Sk\l{}odowska-Curie grant agreement N\textsuperscript{\underline{o}} 945322.  \includegraphics[scale=0.025]{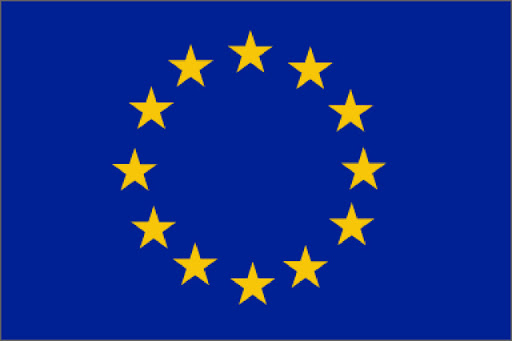}	Moreover, he was partially supported by GeoMod AAPG2019 (ANR-DFG), `Geometric and Combinatorial Configurations in Model Theory'.}

\subjclass[2020]{Primary 03C60, Secondary 03C66, 12L12, 03C45}

\keywords{model theory, difference fields, higher amalgamation, additive character}

\begin{document}

\begin{abstract} 
We show two results on higher amalgamation in the theory $\mathrm{ACFA}^{+}$, the model companion of the theory of difference fields with an additive character (added as a continuous logic predicate) on the fixed field in characteristic 0. On one hand, we show that the non-trivial condition for 3-amalgamation established in a preceding paper is not sufficient for 4-amalgamation. On the other hand, we show that when working over substructures whose $\mathcal{L}_{\sigma}$-reduct is a model of $\mathrm{ACFA}$, $n$-amalgamation holds for all $n\geq 3$.
\end{abstract}

\maketitle




\section{Introduction}
This article is a continuation of the study of the model-theoretic properties of the theory $\mathrm{ACFA}^{+}$ that was begun in \cite{ludwig2025modeltheorydifferencefields}. The theory $\mathrm{ACFA}^{+}$ is the model companion of the theory of difference fields with an additive character (added as a continuous logic predicate) on the fixed field \cite[Theorem 3.19]{ludwig2025modeltheorydifferencefields}. One of the principal motivations to study $\mathrm{ACFA}^{+}$ is that it is the common theory (in characteristic $0$) of the algebraic closure of finite fields $\bar{\mathbb{F}}_{q}$ together with the Frobenius automorphism $\mathrm{Frob}_{q}$ and an additive character $\Psi_{q}$ on $\mathbb{F}_{q}$ \cite[Theorem 3.28]{ludwig2025modeltheorydifferencefields}. In his work on pseudofinite fields with an additive character \cite{Hrushovski2021AxsTW} Hrushovski proposes to study the theory $\mathrm{ACFA}^{+}$ and already observes that the Kim-Pillay group of (completions of) $\mathrm{ACFA}^{+}$ is no longer necessarily totally disconnected as opposed to (what is knwon for) simple theories in classical logic.
However, he conjectures that it is always abelian \cite[Section 6.3]{Hrushovski2021AxsTW}. In \cite[Theorem 5.11]{ludwig2025modeltheorydifferencefields} this abelianity is proved with the main ingredient being the characterisation of those algebraically closed difference subfields over which the Independence Theorem holds in $\mathrm{ACFA}^{+}$. Building on this characterisation, simplicity of $\mathrm{ACFA}^{+}$ is obtained as a straight-forward consequence and moreover a full description of the continuous logic imaginaries present in $\mathrm{ACFA}^{+}$ is deduced \cite[Section 6]{ludwig2025modeltheorydifferencefields}.\\ 
The Independence Theorem can be equivalently formulated in terms of 3-amalgamation of types. This notion has higher-dimensional analogues, referred to as \textit{n-amalgamation}, which date back to Shelah  \cite{shelahclassificationtheory}. An $n$-amalgamation problem (Definition \ref{definitionnamalgamationproblem}) can be pictured as a simplicial complex of types. While uniqueness of 2-amalgamation in stable theories (stationarity) and 3-existence in simple theories (the Independence theorem) are cornerstones of model-theoretic classification theory, the higher dimensional analogues are less well developed and prominent. However, an early application of 4-amalgamation is in the context of the group configuration theorem in simple theories \cite{constructinghyperdefgroup}. In a stable theory, 4-amalgamation is linked to definable groupoids, similar to how 3-amalgamation is linked to imaginaries. This was established in \cite{grupoidshrushovski}, leading to subsequent generalisations, such as, for example \cite{typeamalgamationandpolygroupoids}, \cite{Goodrick2010AmalgamationFA}
 and \cite{grupoidscoversgoodrickkolesnikov}. An account of the homological flavour of higher amalgamation is given in \cite{homologygroupsoftypesgoodrickkimkolesnikov}. An interesting but often intricate question is to determine whether non-trivial higher-dimensional phenomena in a given theory appear or if n-amalgamation already holds whenever all obstructions to 3-amalgamation are eliminated. While the latter holds for many natural theories of (enriched) fields such as $\mathrm{ACF}$ (algebraically closed fields), $\mathrm{DCF}$ (differentially closed fields) or $\mathrm{ACFA}$ (existentially closed difference fields), a natural example for the first case arises from the context of compact complex manifolds \cite{complexcompactmanifoldswithautom} where 4-amalgamation does not hold over general algebraically closed (including imaginaries) sets.\\
As announced in \cite[7.1]{ludwig2025modeltheorydifferencefields} higher amalgamation in $\mathrm{ACFA}^{+}$ turns out to be surprisingly intricate as compared to $\mathrm{ACFA}$ where the results on 3-amalgamation generalise in a straight-forward manner to higher dimensions. This article will be dedicated to the proof of two results in opposite directions: On one hand, we show that the condition on 3-amalgamation obtained in \cite{ludwig2025modeltheorydifferencefields} is already not sufficient to control 4-amalgamation (Section \ref{sectioncounterexample}). On the other hand, when working over models, $n$-amalgamation holds for all $n\in\mathbb{N}$ due to a stability-theoretic argument (Section \ref{sectionnamalgovermodelschapterfour}). In particular, we can not simply \textit{lift} a counterexample to 3-amalgamation to a counterexample to 4-amalgamation over a model.

\subsection*{Set-up and previous results} 
We recall from \cite{ludwig2025modeltheorydifferencefields} that the language $\mathcal{L}_{\sigma,+}$ consists of the language of rings $\mathcal{L}_{\mathrm{ring}}$ together with a unary function symbol $\sigma$ and an $S^{1}\cup\{0\}$-valued continuous logic predicate $\Psi$. As stated above $\mathrm{ACFA}^{+}$ is the model companion of the theory of $\mathcal{L}_{\sigma,+}$-structures given by difference fields with an additive character on the fixed field \cite[Theorem 3.19]{ludwig2025modeltheorydifferencefields}. Independence in $\mathrm{ACFA}^{+}$ is directly inherited from its reduct $\mathrm{ACFA}$ (see Definition \ref{definitionindependenceacfa}) and is given by algebraic independence (after passing to the $\sigma$-algebraic closure $\mathrm{acl}_{\sigma}$ of the relevant sets). The main results from \cite{ludwig2025modeltheorydifferencefields} that we will use throughout this paper are Lemma \ref{lemmatorsornotrealisedinalgclosure} (which we explicitly present below after Theorem \ref{theoremintroA}) and the following description of 3-amalgamation:\\
In $\mathrm{ACFA}^{+}$ 3-amalgamation holds over some $A=\mathrm{acl}_{\sigma}(A)$ if and only if $A$ is $\sigma$-AS-closed, that is, for all $a\in A$, there is $b\in A$ with $\sigma(b)-b=a$ \cite[Theorem 4.18]{ludwig2025modeltheorydifferencefields}. So, while the presence of the continuous logic predicate $\Psi$ is the reason for this non-trivial condition for 3-amalgamation, the condition itself is purely formulated in the language $\mathcal{L}_{\sigma}$ of difference fields. This construction together with some technical general results on higher amalgamation are recalled in more detail in Section \ref{sectionpreliminaries}.

\subsection*{Presentation of results}
In Section \ref{sectioncriterionhigheramalg} we translate higher amalgamation in $\mathrm{ACFA}^{+}$ into a condition that is purely formulated in the language $\mathcal{L}_{\sigma}$ of difference fields (Lemma \ref{lemmacorrespondingvectorpsaceproperty}). From there on, the continuous logic appears only implicitly and the main technical results concern in a certain sense purely the model theory of difference fields.\\
The first main result deals with higher amalgamation over a model. Given the non-trivial condition on 3-amalgamation, one could wonder whether 4-amalgamation over models has to fail simply by constructing a 4-amalgamation problem where (the realisation of) one corner $\bar{a}_{1}$ is not $\sigma$-AS-closed and the types for the other corners are copies of a 3-amalgamation problem that does not have a solution over $\bar{a}_{1}$. This argument however does not work as the independence condition would be violated and the following result even implies that no such counterexample can be constructed.
\begin{theoremintro}(See Theorem \ref{theoremnamalgamationovermodels}.)\label{theoremintroA}
In $\mathrm{ACFA}^{+}$ $n$-amalgamation holds over all substructures whose $\mathcal{L}_{\sigma}$-reduct is a model of $\mathrm{ACFA}$.
\end{theoremintro}

The idea of the proof is to first generalise the notion of $\sigma$-AS-closedness to a higher dimensional condition (Definition \ref{definitionntorsorclosed}) which implies $n$-amalgamation and then to show using a stability-theoretic argument that any $A\models \mathrm{ACFA}$ satisfies this condition for every $n\in\mathbb{N}$. The heart of the argument is the following key lemma, which makes the stability-theoretic arguments accessible.

\begin{lemma*}
Let $(K,\sigma)$ be a difference field of characteristic $0$ and $a\in K$. Assume that there is no $b\in K$ with $\sigma(b)-b=a$, then there is also no $b\in K^{\mathrm{alg}}$ as such.    
\end{lemma*}

It has a completely elementary proof and was already present in the author's previous work on $\mathrm{ACFA}^{+}$ (see \cite[Lemma 4.12]{ludwig2025modeltheorydifferencefields}). Interestingly, two variants of it (Lemma \ref{lemmalineardiffequationnotrealisedinalgclosure} and \ref{lemmamultiplicativediffequationnotrealisedinalgclosure}) are also used in the subsequent Section \ref{sectioncounterexample}.\\
Considering the result on higher amalgamation over models, it is natural to wonder whether being $\sigma$-AS-closed is sufficient to ensure $n$-amalgamation for $n>3$ as well. In particular, given that moreover in $\mathrm{ACFA}$, $\mathrm{PF}^{+}$ and many more natural theories $n$-amalgamation behaves (and is proved) very similarly to 3-amalgamation, it would be natural to assume that $\sigma$-AS-closed indeed is sufficient. The second main results of this article states that this is not true and that already for 4-amalgamation $\sigma$-AS-closedness is not a sufficient condition. More precisely, we show the following

\begin{theoremintro}(See Theorem \ref{theoremnot4amalagamation}.)
There is some model of $\mathrm{ACFA}^{+}$ containing a $\sigma$-AS-closed set over which 4-amalgamation does not hold.
\end{theoremintro}
The proof has two main components. First, one shows that it suffices to construct a $\sigma$-AS-closed set $E$ that is not 2-$\sigma$-AS-closed (see Definition \ref{definitionntorsorclosed}). Next, using the aforementioned generalisations (Lemma \ref{lemmalineardiffequationnotrealisedinalgclosure} and \ref{lemmamultiplicativediffequationnotrealisedinalgclosure}) of the above lemma one constructs such an $E$ using a chain argument avoiding realisations of certain difference equations (see Lemma \ref{theoremcounterexampleexists}). The main idea behind the construction is a simple calculation of difference equations (see the proof of Theorem \ref{theoremnot4amalagamation}). In contrast to the case of 3-amalgamation, not only Torsors of the additive group are present, but multiplicative twists (solutions to equations of the form $\sigma(x)-ex=e$) play an essential role.


\subsection*{Structure of the article}Section \ref{sectionpreliminaries} is mostly devoted to recall necessary results from \cite{ludwig2025modeltheorydifferencefields} on the theory $\mathrm{ACFA}^{+}$ as well as on general higher amalgamation. In Section \ref{sectioncriterionhigheramalg} $n$-amalgamation is translated into a condition in the language of difference fields which we then use in the following two sections that contain the proofs of the two main results of this article. In Section \ref{sectionnamalgovermodelschapterfour} it is shown that $n$-amalgamation holds over substructures whose $\mathcal{L}_{\sigma}$-reduct is a model of $\mathrm{ACFA}$ and in Section \ref{sectioncounterexample} we prove that $\sigma$-AS-closedness is not sufficient for 4-amalgamation. Finally, we state some remarks and natural open questions in Section \ref{sectionquestionsandremarks}.

\subsection*{Acknowledgments}
This article is based on Chapter 4 of my PhD thesis \cite{ludwig:tel-05236078} written at the École Normale Supérieure Paris. I would like to thank my supervisors Zoé Chatzidakis and Martin Hils for their invaluable support during my thesis. Moreover, I thank an anonymous referee for useful comments and suggestions. 

\section{Preliminaries}\label{sectionpreliminaries}

\subsection{The theory $\mathrm{ACFA}^{+}$}

We now recall some of the main notions and results around the theory $\mathrm{ACFA}^{+}$ as established in \cite{ludwig2025modeltheorydifferencefields}. In particular, the characterisation of 3-amalgamation which we recall in Section \ref{sectionprelim3amalgam} will be indispensable for our results on higher amalgamation. We start by recalling the axiomatisation of the theory $\mathrm{ACFA}$, the model companion of the theory of difference fields. For us a difference field will be a pair $(K,\sigma)$ consisting of a field together with an automorphism $\sigma:K\rightarrow K$. Let $\mathcal{L}_{\sigma}$ be the language that consists of the ring language together with a unary function symbol $\sigma$.

\begin{notation}
    In the following we will always write $\mathrm{acl}(\cdot)$ for the model-theoretic algebraic closure and $H^{\mathrm{alg}}$ or $\bar{H}$ for the field-theoretic algebraic closure of some subfield $H\subseteq K$.
\end{notation}

\begin{definition}(1.1 in \cite{acfa})
    The $\mathcal{L}_{\sigma}$-theory $\mathrm{ACFA}$ is given by a scheme of axioms expressing
the following properties of the $\mathcal{L}_{\sigma}$-structure $(K,\sigma)$:
\begin{itemize}
    \item $\sigma$ is an automorphism of $K$.
    \item $K$ is an algebraically closed field.
    \item For every absolutely irreducible varieties $V,U$ with $V\subseteq U\times \sigma(U)$ projecting generically onto $U$ and $\sigma(U)$ and every proper algebraic subset $W\subset V$ there is $a\in U(K)$ such that $(a,\sigma(a))\in V\backslash W$.
\end{itemize}

\end{definition} 

\begin{fact}\label{basicfactsonacfa}(1.1-1.3,1.11 in \cite{acfa})
\begin{itemize}
    \item Every difference field embeds in a model of $\mathrm{ACFA}$ and moreover $\mathrm{ACFA}$ is model-complete.
    \item For two models $(K_{1},\sigma_{1}),(K_{2},\sigma_{2})$ of $\mathrm{ACFA}$ with a common difference subfield $E$, we have \[(K_{1},\sigma_{1})\equiv_{E}(K_{2},\sigma_{2})\;\iff\; \left(E^{\mathrm{alg}},\sigma_{1}\restriction_{E^{\mathrm{alg}}}\right)\cong\left(E^{\mathrm{alg}},\sigma_{2}\restriction_{E^{\mathrm{alg}}}\right).\]
    \item The fixed field $F:=\mathrm{Fix}(\sigma)$ of any model of $\mathrm{ACFA}^{+}$ is a pseudofinite field, i.e., it is perfect, PAC and $\mathrm{Gal}(F)=\hat{\mathbb{Z}}$.
    \item  $F$ is stably embedded in any model $K\models\mathrm{ACFA}$. The restriction of every $\mathcal{L}_{\sigma}(K)$-definable subset to $F$ is definable with parameters from $F$ only using the language $\mathcal{L}_{\mathrm{ring}}$.
\end{itemize}
    
\end{fact}

Note that the following is natural to assume when we work later with $\mathrm{ACFA}^{+}$ as an additive character on a pseudofinite field in positive characteristic is already interpretable in the pure field structure.

\begin{convention}
From now on we fix to work in characteristic $0$. In particular, $\mathrm{ACFA}$ is assumed to contain a set of sentences stating that the characteristic is $0$.
\end{convention}

\begin{notation}
    We write $S^{1}$ for the unit circle and $\mathbb{T}^{n}\cong S^{1}\times\cdots\times S^{1}$ for the $n$-dimensional complex torus.
\end{notation}

\begin{definition}
    We denote by $\mathcal{L}_{\sigma,+}$ the extension of the language $\mathcal{L}_{\sigma}$ by a unary continuous logic predicate $\Psi$ which is allowed to take values in $S^{1}\cup\{0\}$.
\end{definition}

Note that while we work in continuous logic, we only add a continuous logic predicate to the language. So, in particular, equality is treated in the usual way as in classical discrete logic (or, in other words, the underlying metric of the structure is the discrete metric). As outlined in the introduction of \cite{Hrushovski2021AxsTW} this embeds in the usual presentation of real-valued continuous logic as in \cite{mtfms} by taking two predicates, one for the real and one for the imaginary part of $\Psi$.

\begin{notation}

We write $\Psi^{(n)}:K^{n}\rightarrow\{S^{1}\cup\{0\}\}^{n}$ for the map $(x_{1},\dots,x_{n})\rightarrow(\Psi(x_{1}),\dots,\Psi(x_{n}))$.
\end{notation}

\begin{definition}
   A rational hyperplane over $F$ (in $\mathbb{A}^{n}$) is a variety that is defined by an equation of the form $\sum_{1\leq i\leq n}z_{i}X_{i}=b$ where $b\in F$ and $z_{i}\in\mathbb{Z}$ for all $1\leq i\leq n$. Moreover, we require that $z_{i}\neq 0$ for some $1\leq i\leq n$. The rational hyperplane is said to have height $\leq m$, if $|z_{i}|\leq m$ for all $1\leq i\leq n$.
\end{definition}

We now state  \cite[Definition 3.4]{ludwig2025modeltheorydifferencefields}. Note that it was already proposed by Hrushovski in \cite[Section 6.3]{Hrushovski2021AxsTW} to investigate the theory $\mathrm{ACFA}^{+}$. We now state its definition which, as we explain afterwards, directly emerges from the theory $\mathrm{PF}^{+}$ from \cite{Hrushovski2021AxsTW}.

\begin{definition}(Definition 3.4 in \cite{ludwig2025modeltheorydifferencefields})\label{definitionacfaplus}
The theory $\mathrm{ACFA}^{-}$ consists of the $\mathcal{L}_{\sigma}$-theory $\mathrm{ACFA}$ together with axioms stating that $\Psi(K\backslash F)=0$ for $F:=\mathrm{Fix}(\sigma)$ as well as that $\Psi\restriction_{F}:(F,+)\rightarrow(S^{1},\cdot)$ is a group homomorphism. The theory $\mathrm{ACFA}^{+}$ then extends the above theory $\mathrm{ACFA}^{-}$ by the following set of axioms:\\

\noindent
($\star$) Let $n,m\in\mathbb{N}$. Let $h\in\mathbb{Q}[z_{1},z_{1}^{-1},\dots,z_{n},z_{n}^{-1}]$ be a finite Fourier series (Laurent polynomial) with degree (in every $z_{i}$) bounded by $\leq m$ which is moreover real-valued on $\mathbb{T}^{n}$ and does not have a constant term. For any absolutely irreducible curve $C$ that is defined over $F$ such that $C\subset\mathbb{A}^{n}$ is not contained in any rational hyperplane over $F$ of height at most $m$, the following holds
\[\sup\{h(\Psi^{(n)}(\bar{x}))\;:\;\bar{x}\in C(F)\}\geq 0.\]

\end{definition}

\begin{fact}\label{factaxiomisexpressible}(Lemma 3.5 in \cite{Hrushovski2021AxsTW}.)\label{factequivalenceaxiomtodensevaluesoncruve}
   Let $F$ be a field of characteristic $0$ and $\Psi$ an additive character on $F$. Then, the condition $(\star)$ from Definition \ref{definitionacfaplus} holds in $F$ if and only if the following is true:\\
   Given any absolutely irreducible curve $C\subset \mathbb{A}^{n}$ which is defined over $F$ and is not contained in a rational hyperplane over $F$, the set $\Psi^{(n)}(C(F))$ is a dense subset of $\mathbb{T}^{n}$ (in the euclidean topology).
   
\end{fact}

The above axiom $(\star)$ together with the conditions that $\mathrm{Gal}(F)=\hat{\mathbb{Z}}$ and $\Psi$ being an additive character on $F$ (as well as $\mathrm{char}(F)=0$) characterise the models of the theory $\mathrm{PF}^{+}$ that was introduced by Hrushovski in \cite{Hrushovski2021AxsTW}. Hrushovski's paper is crucial for the prior work on $\mathrm{ACFA}^{+}$ in \cite{ludwig2025modeltheorydifferencefields}, among other things he shows that $\mathrm{PF}^{+}$ is the common theory (in characteristic $0$) of finite fields with non-trivial additive character as well as the definability of the Chatzidakis-van den Dries-Macintyre counting measure in $\mathrm{PF}^{+}$. Thus, in particular, $\mathrm{ACFA}^{+}$ amalgamates the theory $\mathrm{ACFA}$ and the theory $\mathrm{PF}^{+}$ from \cite{Hrushovski2021AxsTW} in the following sense. An $\mathcal{L}_{\sigma,+}$-structure is a model of $\mathrm{ACFA}^{+}$ if and only if its $\mathcal{L}_{\sigma}$-reduct is a model of $\mathrm{ACFA}$ and the $\mathcal{L}_{+}=\mathcal{L}_{\mathrm{ring}}\cup\{\Psi\}$-structure on the fixed field is a model of $\mathrm{PF}^{+}$.\\
One of the motivations to study $\mathrm{ACFA}^{+}$ is that the above mentioned result of Hrushovski for $\mathrm{PF}^{+}$ in \cite{Hrushovski2021AxsTW} and the corresponding result for $\mathrm{ACFA}$ in \cite{elementaryfrobenius} can be combined to obtain the following.

\begin{fact}(Theorem 3.28 in \cite{ludwig2025modeltheorydifferencefields})
The $\mathcal{L}_{\sigma,+}$-theory $\mathrm{ACFA}^{+}$ is the common theory $\mathrm{Th}((\bar{\mathbb{F}}_{q},\mathrm{Frob}_{q},\Psi_{q})_{q})$ of the algebraic closure of finite fields $\bar{\mathbb{F}}_{q}$ with the Frobenius $\mathrm{Frob}_{q}:\bar{\mathbb{F}}_{q}\rightarrow\bar{\mathbb{F}}_{q}$ and a non-trivial additive character $\Psi_{q}$ on $\mathbb{F}_{q}$ together with a set of axioms stating that the characteristic is $0$. (Here $q$ ranges over all prime powers.)
    
\end{fact}

The following lemma from \cite{ludwig2025modeltheorydifferencefields} will be of use later on in our results on higher amalgamation.

\begin{fact}(Lemma 4.4 in \cite{ludwig2025modeltheorydifferencefields})\label{factlinindconsistent}
  Let $K\models \mathrm{ACFA}^{+}$ and $E=\mathrm{acl}_{\sigma}(E)\subseteq K$. Let $\bar{c},\bar{d}$ be tuples in $K\backslash E$ and, moreover, $\bar{c}\in F^{n}$. Assume that $\bar{c}$ is not contained in any rational hyperplane over $F\cap \mathrm{acl}_{\sigma}(E\bar{d})$. Let $p(\bar{x},\bar{y}):=\mathrm{tp}_{\mathcal{L}_{\sigma}}(\bar{c}\bar{d}/E)\cup \mathrm{tp}(\bar{d}/E)$.
 For any $\bar{k}\in\mathbb{N}^{n}$ and any $\bar{r}\in\mathbb{T}^{n}$, the partial $\mathcal{L}_{\sigma,+}$-type $p(\bar{x},\bar{y})$ is consistent with $\Psi^{(n)}(\frac{1}{k_{1}}x_{1},\dots,\frac{1}{k_{n}}x_{n})=\bar{r}$.  
\end{fact}

\subsection{Independent type-amalgamation}\label{sectionindtype-amalgfirstchapter}
In this section we present independent type-amalgamation, that is, what we refer to as \textit{higher amalgamation}. We present some variants that will prove useful later on. The main source is \cite{grupoidshrushovski}. At the end of the section we give the technical Lemma \ref{lemmadecomposabilityofgeneraladditiveequations} to which we will later often refer as \textit{decomposability of additive equations}. While it follows from basic considerations of higher amalgamation in $\mathrm{ACF}$, it will prove very useful when we will investigate higher amalgamation in the theory $\mathrm{ACFA}^{+}$.

\begin{convention}
   In this section, we fix to work in the following setting: $\mathcal{L}$ is a (possibly continuous) expansion of the language of rings $\mathcal{L}_{ring}$. We fix $T$ to be a complete $\mathcal{L}$-theory extending the theory $\mathrm{ACF}$ of algebraically closed fields. We assume that $T$ carries a notion of independence that implies independence in $\mathrm{ACF}$ and further satisfies all the axioms of the Kim-Pillay characterisation of simple theories via an independence relation (see Theorem 1.51 in \cite{simplicityCATs} for the statement in the context of CATs which comprises continuous logic) possibly except of the Independence theorem (axiom (8) in Theorem 1.51 in \cite{simplicityCATs}). Moreover, we fix $\mathcal{M}$ to be a sufficiently saturated model of $T$ and $E\subseteq M$ to be a small algebraically closed set.
\end{convention}

\begin{remark}
    For most notions in this section, it would suffice to work with an arbitrary theory that carries a notion of independence satisfying axioms (1)-(7) in Theorem 1.51 in \cite{simplicityCATs}. The underlying field structure does not make an appearance before Lemma \ref{lemmaspecialisationringformulas}.\footnote{And even there it would suffice to work with an underlying stable theory satisfying $n$-uniqueness (see Chapter 4 of \cite{grupoidshrushovski} for a Definition) for every $n\in\mathbb{N}$ instead of $\mathrm{ACF}$.} 
\end{remark}

\begin{notation}
We denote by $\mathbf{n}$ the set $\{1,\dots,n\}$ and write $\mathcal{P}^{-}(\mathbf{n})=\mathcal{P}(\mathbf{n})\backslash\{\mathbf{n}\}$. Further, for all $w\in\mathcal{P}^{-}(\mathbf{n})$ let $\bar{x}_{w}$ be a (possibly infinite) tuple of variables such that $\bar{x}_{w}\cap \bar{x}_{w^{\prime}}=\bar{x}_{w\cap w^{\prime}}$ for all $w^{\prime}\in \mathcal{P}^{-}(\mathbf{n})$.

\end{notation}
\begin{definition}\label{definitionnamalgamationproblem}
    We consider a system $(p_{w})_{w\in\mathcal{P}^{-}(\mathbf{n})}$, where $p_{w}(\bar{x}_{w})$ is a complete $T$-type over $E$. We call $(p_{w})_{w\in\mathcal{P}^{-}(\mathbf{n})}$ an \textit{n-amalgamation problem (over $E$)} if the following conditions hold:
    \begin{enumerate}
        \item For all $w,w^{\prime}\in \mathcal{P}^{-}(\mathbf{n})$, if $w\subseteq w^{\prime}$, then $p_{w}(\bar{x}_{w})\subseteq p_{w^{\prime}}(\bar{x}_{w^{\prime}})$.
        \item For any $\Tilde{w}\in \mathcal{P}^{-}(\mathbf{n})$, for any $w,w^{\prime}\subseteq\Tilde{w}$, for any $\bar{a}_{\Tilde{w}}\models p_{\Tilde{w}}$ we have that $\bar{a}_{w}$ is independent from $\bar{a}_{w^{\prime}}$ over $\bar{a}_{w\cap w^{\prime}}=\bar{a}_{w}\cap \bar{a}_{w^{\prime}}$ where $\bar{a}_{w}\subseteq \bar{a}_{\Tilde{w}}$ is the subtuple corresponding to $\bar{x}_{w}\subseteq \bar{x}_{\Tilde{w}}$ (and similarly for $w^{\prime},w\cap w^{\prime}$). 
        \item If $\bar{a}_{w}\models p_{w}$, then $\bar{a}_{w}$ enumerates $\mathrm{acl}(E,\bar{a}_{\{i\}}\,|\,i\in w)$.
        
    \end{enumerate}
    We call a complete type $p_{\mathbf{n}}(\bar{x}_{\mathbf{n}})$ a \textit{solution} to the $n$-amalgamation problem $(p_{w})_{w\in\mathcal{P}^{-}(\mathbf{n})}$, if the system $(p_{w})_{w\in\mathcal{P}(\mathbf{n})}$ still fulfills the above conditions 1.-3.
    We say that the theory $T$ has \textit{n-amalgamation} over a class $\mathcal{C}$ of algebraically closed sets if for every $m\leq n$ every $m$-amalgamation problem over some $E\in\mathcal{C}$ (with $E\subseteq\mathcal{M}\models T$) has a solution. An incomplete theory has $n$-amalgamation if all of its completions have it.
\end{definition}

\begin{remark}
    When we will later consider higher amalgamation in the incomplete theory $\mathrm{ACFA}^{+}$ we will always work in an arbitrary but fixed completion of $\mathrm{ACFA}^{+}$.
\end{remark}

\begin{definition}
A type $p_{w}$ from a system $(p_{w})_{w\in\mathcal{P}^{-}(\mathbf{n})}$ will be called an $r$-type, if $|w|=r$.
\end{definition}

\begin{notation}
    For $i\in w$ we write $\bar{x}_{i}$ (or $\bar{a}_{i}$) instead of $\bar{x}_{\{i\}}$ (or $\bar{a}_{\{i\}}$).
    Let $n\in\mathbb{N}$. We set $\widehat{w}_{{i}}:=\{1,\dots,n\}\backslash\{i\}$ for $1\leq i\leq n$. Moreover, we will write $\widehat{w}_{i,j}$ for $\{1,\dots,n\}\backslash\{i,j\}$ and $\widehat{w}_{J}$ for $\{1,\dots,n\}\backslash J$ for some $J\subseteq \{1,\dots,n\}$. As $w$ (resp. $\widehat{w}_{i}$) will be reserved for indices of elements of a given system (of types, tuples,...) we will write $\mathcal{P}(\widehat{\mathbf{n}}_{i})$ for $\mathcal{P}(\widehat{w}_{i})$.
\end{notation}

\begin{definition}
    We call a system of tuples $(\bar{a}_{w})_{w\in\mathcal{P}(\mathbf{n})}$ an \textit{independent n-system over $E$} if we have for any $w,w^{\prime}\in \mathcal{P}(\mathbf{n})$ that $\bar{a}_{w}$ is independent from $ \bar{a}_{w^{\prime}}$ over $\bar{a}_{w\cap w^{\prime}}$ and, moreover, that $\bar{a}_{w}$ enumerates $\mathrm{acl}(E,\bar{a}_{i}\,|\,i\in w)$. We always assume $E\subsetneq \bar{a}_{i}$ for all $1\leq i\leq n$.
\end{definition}

\begin{remark}
    
Any solution of an $n$-amalgamation problem over $E$ yields an independent $n$-system over $E$.
\end{remark}
We can extend Definition \ref{definitionnamalgamationproblem} to a \textit{fragmentary amalgamation problem} in the following way:
\begin{definition}\label{definitionn-m-amalgam} 
    Let $m<n$ and denote by $\mathcal{P}_{\leq m}(\mathbf{n})$ the set of $X\subseteq\{1,\dots,n\}$ with $|X|\leq m$. Then we call a system of complete $T$-types $(p_{w})_{w\in \mathcal{P}_{\leq m}(\mathbf{n})}$ an $(m,n)$-amalgamation problem if it satisfies the same conditions as in Definition \ref{definitionnamalgamationproblem} and we say that it has a solution if it can be completed to a system $(p_{w})_{w\in\mathcal{P}(\mathbf{n})}$ still satisfying the conditions. As before we say that the theory $T$ \textit{has (m,n)-amalgamation over a class $\mathcal{C}$ of algebraically closed sets} if every $(m,n)$-amalgamation problem over some $E\in\mathcal{C}$ has a solution.
\end{definition}

\begin{remark}
    Note that an $(n-1,n)$-amalgamation problem is simply an $n$-amalgamation problem as in Definition \ref{definitionnamalgamationproblem}.
\end{remark}

\begin{fact}\label{factamalgamationimpliespartiaamalgamation}(Lemma 4.1 (2) in \cite{grupoidshrushovski}.)
    If $T$ has $n$-amalgamation over all algebraically closed sets, then it also has $(n-1,n+k)$-amalgamation for all $k\geq 1$ over all algebraically closed sets.
\end{fact}
Next, we have to introduce a further technical notion, that of a \textit{partial} amalgamation problem. It will be used to state in Lemma \ref{lemmapartialamalgamation} the simple fact that once we work in a theory with independent $n$-amalgamation, then we can also amalgamate if the types are only determined for a subset of $\mathcal{P}_{\leq n-1}(\mathbf{n+k+1})$. 
\begin{definition}\label{definitionpartialamalgamation}
    Let $\emptyset\neq Y\subsetneq \mathcal{P}_{\leq m}(\mathbf{n})$ be downward-closed, i.e., if $X_{1}\subseteq X_{2}\in Y$, then $X_{1}\in Y$ and further containing all one-element subsets. We call a system of types $(p_{w})_{w\in Y}$ (over $E$) a \textit{partial} $(m,n)$-amalgamation problem (over $E$) if it satisfies the conditions of Definition \ref{definitionnamalgamationproblem} (quantifying over $Y$ instead of $\mathcal{P}^{-}(\mathbf{n})$).
\end{definition}

\begin{lemma}\label{lemmapartialamalgamation}
    If $T$ has $ n$-amalgamation for all algebraically closed sets, then it also has partial $(n-1,n+k)$-amalgamation for all $k\geq 1$ over all algebraically closed sets.
\end{lemma}
\begin{proof}
    Let $m$ be minimal such that $Y\subseteq \mathcal{P}_{\leq m}(\mathbf{n+k})$. It suffices to complete $(p_{w})_{w\in Y}$ to an $(m,n+k)$-amalgamation problem. This can be done iteratively for $2\leq i\leq m$, using $i$-amalgamation in each step, starting with $i=2$:
    In the step for $i=r$, every $W\in \mathcal{P}_{=r}(\mathbf{n+k})\backslash\{w\in Y\,|\, |w|=r\}$ yields an $i$-amalgamation problem $(p_{w})_{w\in W}$  that can be amalgamated to some type $p_{W}$.
\end{proof}

\begin{notation}
We denote by $ p\restriction_{\mathcal{L}_{\mathrm{ring}}}$ the restrictions of a type $p$ in $\mathcal{L}$ to the language $\mathcal{L}_{\mathrm{ring}}$.
\end{notation}

\begin{definition}
We work in an independent $n$-system $(\bar{a}_{w})_{w\in\mathcal{P}(\mathbf{n})}$ over $E$. Let $J_{0},\dots,J_{k}\subset \{1,\dots,n\}$ be sets and $\bar{b}_{J_{0}},\dots,\bar{b}_{J_{k}}$ be tuples such that $\bar{b}_{J_{i}}\in \bar{a}_{J_{i}}$ for every $0\leq i\leq k$. Let $\phi(\bar{x}_{0},\dots,\bar{x}_{k})$ be an $\mathcal{L}$-formula such that $\phi(\bar{b}_{J_{0}},\dots,\bar{b}_{J_{k}})$ holds in $\bar{a}_{\mathbf{n}}$. Then we say that we can \textit{specialise} $\phi(\bar{b}_{J_{0}},\dots,\bar{b}_{J_{k}})$ (to $\bar{a}_{J_{0}}$), if we can find tuples $\bar{b}_{J_{1}\rightarrow J_{0}},\dots,\bar{b}_{J_{k}\rightarrow J_{0}}$ with $\bar{b}_{J_{i}\rightarrow J_{0}}\in \bar{a}_{J_{i}\cap J_{0}}$ for any $1\leq i\leq k$ such that $\phi(\bar{b}_{J_{0}},\bar{b}_{{J}_{1}\rightarrow J_{0}},\dots,\bar{b}_{\Tilde{J}_{k}\rightarrow J_{0}})$ holds.\\
If this is true for any choice of independent $n$-system $(\bar{a}_{w})_{w\in\mathcal{P}(\mathbf{n})}$ and of tuples $\bar{b}_{J_{0}},\dots,\bar{b}_{J_{k}}$ we simply say that $\phi(\bar{x}_{0},\dots,\bar{x}_{k})$ specialises in the variable tuple $\bar{x}_{0}$ over $E$.
\end{definition}
The following is probably standard. We essentially use the same proof as for independent $n$-amalgamation in $\mathrm{ACFA}$ given in (1.9) of \cite{acfa}.
\begin{lemma}\label{lemmaspecialisationringformulas}
     (Recall that we still work in the $\mathcal{L}$-theory $T$.) Every $\mathcal{L}_{\mathrm{ring}}$ formula specialises over $E$.
\end{lemma}
\begin{proof}
    Let $\phi(\bar{x}_{0},\dots,\bar{x}_{k})$ be an $\mathcal{L}_{\mathrm{ring}}$-formula and fix an independent $n$-system $(\bar{a}_{w})_{w\in\mathcal{P}(\mathbf{n})}$ over $E$. Assume that $\phi(\bar{b}_{J_{0}},\dots,\bar{b}_{J_{k}})$ holds for tuples $\bar{b}_{J_{i}}\in \bar{a}_{J_{i}}$ and $J_{0},\dots,J_{k}\subset \{1,\dots,n\}$.\\
    Note that in general for $J,H\subset\{1,\dots,n\}$ with $J\cap H=\emptyset$ we can depict a tuple $\bar{\beta}=\beta_{1},\dots,\beta_{r}$ of elements of the field amalgam of $\bar{a}_{J}$ and $\bar{a}_{H}$ in the following way. We can find finite tuples $\bar{\alpha}_{J}\in \bar{a}_{J}$ and $\bar{\alpha}_{H}\in \bar{a}_{H}$ such that for every $1\leq s\leq r$ we have $\beta_{s}=h_{\beta_{s}}(\bar{\alpha}_{J},\bar{\alpha}_{H})=\frac{f_{\beta_{s}}(\bar{\alpha}_{J},\bar{\alpha}_{H})}{g_{\beta_{s}}(\bar{\alpha}_{J},\bar{\alpha}_{H})}$ for $\mathcal{L}_{\mathrm{ring}}$-terms $f_{\beta_{s}},g_{\beta_{s}}$.\\
    Now instead of general $J,H$, we consider for any $1\leq i\leq k$ the sets $J_{0}\cap J_{i}$ and $J_{i}\backslash J_{0}$. Let $Q_{\bar{b}_{J_{i}}}^{\gamma}(X,\bar{Z})\in\mathbb{Z}[X,\bar{Z}]$ and $\bar{\beta}$ be such that $Q_{\bar{b}_{J_{i}}}^{\gamma}(X,\bar{\beta})$ is the minimal polynomial of an element $\gamma$ from the tuple $\bar{b}_{J_{i}}$ over the field amalgam of $\bar{a}_{J_{i}\cap J_{0}}$ and $\bar{a}_{J_{i}\backslash J_{0}}$. Next, consider the $\mathcal{L}_{\mathrm{ring}}(\bar{a}_{J_{i}\cap J_{0}})$-formula $q_{\bar{b}_{J_{i}}}(\bar{y}_{i},\bar{z})$ which expresses \[\bigwedge_{\gamma\in\bar{b}_{J_{i}}}\;Q_{\bar{b}_{J_{i}}}^{\gamma}\left(y_{i,\gamma},h_{\beta_{1}}(\bar{\alpha}_{J_{0}\cap J_{i}},\bar{z}),\dots,h_{\beta_{r}}(\bar{\alpha}_{J_{0}\cap J_{i}},\bar{z})\right)=0.\]
    Now we return to the initial formula $\phi(\bar{x}_{0},\dots,\bar{x}_{k})$ and find some tuples $\bar{\alpha}_{J_{i}\backslash J_{0}}\in \bar{a}_{J_{i}\backslash J_{0}}$ for all $1\leq i\leq k$ such that \[\mathcal{M}\models\exists(\bar{y}_{i})_{1\leq i\leq k}\;\phi(\bar{b}_{J_{0}},\bar{y})\land\bigwedge_{1\leq i\leq k}q_{\bar{b}_{J_{i}}}(\bar{y}_{i},\bar{\alpha}_{J_{i}\backslash J_{0}}).
    \]
    We use that every $q_{\bar{b}_{J_{i}}}$ is an $\mathcal{L}_{\mathrm{ring}}(\bar{a}_{J_{i}\cap J_{0}})$-formula. Since the type $\mathrm{tp}(\bar{a}_{J_{0}}/\bar{a}_{\widehat{w}_{J_{0}}})\restriction_{\mathcal{L}_{\mathrm{ring}}}$ is the co-heir of $\mathrm{tp}(\bar{a}_{J_{0}}/E)\restriction_{\mathcal{L}_{\mathrm{ring}}}$ we find some tuple $\bar{e}$ from $E$ such that    \[\mathcal{M}\models\exists(\bar{y}_{i})_{1\leq i\leq k}\;\phi(\bar{b}_{J_{0}},\bar{y})\land\bigwedge_{1\leq i\leq k}q_{\bar{b}_{J_{i}}}(\bar{y}_{i},\bar{e}).
    \]
    We define the tuple $(\bar{b}_{{J}_{1}\rightarrow J_{0}},\dots,\bar{b}_{{J}_{k}\rightarrow J_{0}})$ to be a witness of this sentence, i.e., \[\mathcal{M}\models\phi(\bar{b}_{J_{0}},\bar{b}_{{J}_{1}\rightarrow J_{0}},\dots,\bar{b}_{{J}_{k}\rightarrow J_{0}})\land\bigwedge_{1\leq i\leq k}q_{\bar{b}_{J_{i}}}(\bar{b}_{{J}_{i}\rightarrow J_{0}},\bar{e}).
    \]and hence $\bar{b}_{{J}_{i}\rightarrow J_{0}}\in \bar{a}_{J_{i}\cap J_{0}}$ for every $1\leq i\leq k$ which completes the proof.
\end{proof}

We will later have to deal with certain sums that arise in $n$-amalgamation problems. The following is a technical tool that will facilitate the proofs of (higher) amalgamation in the theory $\mathrm{ACFA}^{+}$ later.

\begin{definition}\label{definitionadditiveeauqtion}
We call by an \textit{additive equation of height $n$ over $E$} an equation of the form $\sum_{1\leq i\leq n}b_{\widehat{w}_{i}}=0$ where $b_{\widehat{w}_{i}}\in \bar{a}_{\widehat{w}_{i}}$ for $1\leq i\leq n$ and $(\bar{a}_{w})_{w\in\mathcal{P}(\mathbf{n})}$ is an independent $n$-system over $E$.
\end{definition}

\begin{definition}\label{definitiondecomposableequation}
    Let $\sum_{1\leq i\leq n}b_{\widehat{w}_{i}}=0$ be an additive equation of height $n\geq 3$ over $E$ as in Definition \ref{definitionadditiveeauqtion} with corresponding independent system $(\bar{a}_{w})_{w\in\mathcal{P}(\mathbf{n})}$ and $b_{\widehat{w}_{i}}\in \bar{a}_{\widehat{w}_{i}}$. We say that the equation is \textit{decomposable} if for every $1\leq i\leq n$ and $1\leq j\leq n,\,j\neq i$, there are elements $c_{\widehat{w}_{i}}^{j}\in \bar{a}_{\widehat{w}_{i,j}}$ such that $b_{\widehat{w}_{i}}=\sum_{j\neq i}c_{\widehat{w}_{i}}^{j}$ and $c_{\widehat{w}_{i}}^{j}=-c_{\widehat{w}_{j}}^{i}$ for every pair $i,j\in\{1,\dots,n\}$ where $j\neq i$.
\end{definition}

\begin{notation}\label{remarkindependentsystemonesteplower}
    Let $(\bar{a}_{w})_{w\in\mathcal{P}(\mathbf{n})}$ be an independent $n$-system over $E$. We set for every $w\in \mathcal{P}(\mathbf{n-1})$ the tuple $\Tilde{a}_{w}:=\bar{a}_{w\cup\{n\}}$. Then $(\Tilde{a}_{w})_{w\in\mathcal{P}(\mathbf{n-1})}$ is an independent $(n-1)$-system over the algebraically closed set $\bar{a}_{n}\supset E$. If it is clear from the context what is meant, we will sometimes simply say that $(\bar{a}_{w})_{w\in\mathcal{P}(\mathbf{n-1})}$ is an independent $(n-1)$-system over $a_{n}$. Similarly, we denote by $(\Tilde{a}_{w})_{w\in\mathcal{P}(\widehat{\mathbf{n}}_{i})}$ the corresponding $(n-1)$-system over $\bar{a}_{i}$.
\end{notation}


\begin{lemma}\label{lemmadecomposabilityofgeneraladditiveequations}
Every additive equation of height $n\geq 3$ over $E$ is decomposable.
\end{lemma}

\begin{proof}
Let $\sum_{1\leq i\leq n}b_{\widehat{w}_{i}}=0$ be an additive equation of height $n\geq 3$ over $E$ as in Definition \ref{definitionadditiveeauqtion} with corresponding independent n-system $(\bar{a}_{w})_{w\in\mathcal{P}(\mathbf{n})}$ and $b_{\widehat{w}_{i}}\in \bar{a}_{\widehat{w}_{i}}$.
We will give the proof in two steps: In the first step, we show that for every $1\leq i\leq n$ there are $d_{\widehat{w}_{i}}^{j}$ for every $1\leq j\leq n,\,j\neq i$ such that $b_{\widehat{w}_{i}}=\sum_{j\neq i}d_{\widehat{w}_{i}}^{j}$ and $d_{\widehat{w}_{i}}^{j}\in \bar{a}_{\widehat{w}_{i,j}}$ (but not necessarily $d_{\widehat{w}_{i}}^{j}=-d_{\widehat{w}_{j}}^{i}$).\\
In the second step, we show the statement of the lemma by induction. The main idea will be to apply step 1 to $b_{\widehat{w}_{n}}$ to obtain an additive equation of height $n-1$ over the algebraically closed set $\bar{a}_{n}\supset E$ and then to proceed by induction.\\
\textit{\underline{Step 1:}} 
Consider for every $1\leq i\leq n$ the $\mathcal{L}_{\mathrm{ring}}$-formula \[\phi_{i}(x_{1},\dots,x_{n})\;\;\equiv\;\; x_{i}=-\sum_{1\leq j\leq n,\,j\neq i}x_{j}.\] Then for every $1\leq i\leq n$ the formula $\phi_{i}(b_{\widehat{w}_{1}},\dots,b_{\widehat{w}_{n}})$ holds and we can specialise it to $\bar{a}_{\widehat{w}_{i}}$ by Lemma \ref{lemmaspecialisationringformulas} which yields that there are $d_{\widehat{w}_{i}}^{j}$ such that $b_{\widehat{w}_{i}}=\sum_{j\neq i}d_{\widehat{w}_{i}}^{j}$ and $d_{\widehat{w}_{i}}^{j}\in \bar{a}_{\widehat{w}_{i,j}}$.\\
\textit{\underline{Step 2:}} We start by proving the case $n=3$. Our additive equation then is given by $b_{\widehat{w}_{1}}+b_{\widehat{w}_{2}}+b_{\widehat{w}_{3}}=0$ and by step 1 we obtain that we can write for any triplet $i,k,j$ such that $\{i,k,j\}=\{1,2,3\}$ that $b_{\widehat{w}_{i}}=d_{\widehat{w}_{i}}^{k}+d_{\widehat{w}_{i}}^{j}$ where $d_{\widehat{w}_{i}}^{k}\in a_{j}$ and correspondingly for the other combinations. Now, if we insert this into the additive equation, we obtain for $\delta_{i}:=d_{\widehat{w}_{j}}^{k}+d_{\widehat{w}_{k}}^{j}\in a_{i}$ that $\delta_{1}+\delta_{2}+\delta_{3}=0$ and thus $\delta_{i}\in E$ for all $1\leq i\leq 3$. Then we simply set 
\begin{itemize}
    \item $c_{\widehat{w}_{3}}^{1}=d_{\widehat{w}_{3}}^{1}$
    \item $c_{\widehat{w}_{3}}^{2}=d_{\widehat{w}_{3}}^{2}$
    \item $c_{\widehat{w}_{2}}^{3}=d_{\widehat{w}_{2}}^{3}-\delta_{1}$
    \item $c_{\widehat{w}_{2}}^{1}=d_{\widehat{w}_{2}}^{1}+\delta_{1}$
    \item $c_{\widehat{w}_{1}}^{3}=d_{\widehat{w}_{1}}^{3}-\delta_{2}$
    \item $c_{\widehat{w}_{1}}^{2}=d_{\widehat{w}_{1}}^{2}+\delta_{2}=d_{\widehat{w}_{1}}^{2}-\delta_{1}-\delta_{3}$
\end{itemize}
For this constellation, it follows that $c_{\widehat{w}_{i}}^{j}=-c_{\widehat{w}_{j}}^{i}$ and $b_{\widehat{w}_{i}}=c_{\widehat{w}_{i}}^{j}+c_{\widehat{w}_{i}}^{k}$ for any triplet with $\{i,k,j\}=\{1,2,3\}$. In other words, we have proved the decomposability of $b_{\widehat{w}_{1}}+b_{\widehat{w}_{2}}+b_{\widehat{w}_{3}}$.\\
\textit{Induction step: }We apply Step 1 to $b_{\widehat{w}_{n}}$ and obtain some $d_{\widehat{w}_{n}}^{j}$ for $1\leq j\leq n-1$ such that $b_{\widehat{w}_{n}}=\sum_{1\leq j\leq n-1}d_{\widehat{w}_{n}}^{j}$. Next, for any $1\leq j\leq n-1$, we set $c_{\widehat{w}_{n}}^{j}:=d_{\widehat{w}_{n}}^{j}$ and $e_{\widehat{w}_{j}}:=b_{\widehat{w}_{j}}+c_{\widehat{w}_{n}}^{j}$. Consequently, we have $e_{\widehat{w}_{j}}\in \bar{a}_{\widehat{w}_{j}}$ and from the initial additive equation $\sum_{1\leq i\leq n}b_{\widehat{w}_{i}}=0$ it follows that $\sum_{1\leq j\leq n-1}e_{\widehat{w}_{j}}=0$ holds. But now the latter is an additive equation of height $n-1$ over $\bar{a}_{n}$ and by the induction hypotheses we find for every $1\leq j,k\leq n-1,\,k\neq j$ some $c_{\widehat{w}_{j}}^{k}\in \bar{a}_{w_{j,k}}$ such that $e_{\widehat{w}_{j}}=\sum_{1\leq k\leq n-1,k\neq j}c_{\widehat{w}_{j}}^{k}$ and $c_{\widehat{w}_{j}}^{k}=-c_{\widehat{w}_{k}}^{j}$. We set $c_{\widehat{w}_{j}}^{n}:=-c_{\widehat{w}_{n}}^{j}$ and it follows for any $1\leq j\leq n-1$ that
\[b_{\widehat{w}_{j}}=e_{\widehat{w}_{j}}-c_{\widehat{w}_{j}}^{j}=\sum_{1\leq k\leq n,k\neq j}c_{\widehat{w}_{j}}^{k}\]
Thus, the system $\left(c_{\widehat{w}_{j}}^{k}\right)_{1\leq j,k\leq n,k\neq j}$ yields the decomposabilty of $\sum_{1\leq i\leq n}b_{\widehat{w}_{i}}$, which completes the proof.
\end{proof}

\subsection{3-amalgamation in $\mathrm{ACFA}^{+}$}\label{sectionprelim3amalgam}

We will now state the results on 3-amalgamation in $\mathrm{ACFA}^{+}$ that were obtained in \cite{ludwig2025modeltheorydifferencefields} and will be needed in the following sections.

\begin{notation}
   For a subset $A$ of some model of $\mathrm{ACFA}$ (or $\mathrm{ACFA}^{+}$) we write $\mathrm{acl}_{\sigma}(A)$ to denote the field-theoretic algebraic closure of $\mathrm{cl}_{\sigma}(A)$, the difference field generated by $A$. This is the model-theoretic algebraic closure in $\mathrm{ACFA}$. (See (1.7) in \cite{acfa}.)
    We fix to write $p_{\mathcal{L}_{\sigma}}(\bar{x})$ for $p(\bar{x})\restriction_{\mathcal{L}_{\sigma}}$ given an $\mathcal{L}_{\sigma}^{+}$-type $p(\bar{x})$.
\end{notation}

\begin{definition}\label{definitionindependenceacfa}
  For subsets $A,B,C$ of some model $\mathcal{M}\models \mathrm{ACFA}$ (or $\mathrm{ACFA}^{+}$) such that $C\subseteq A,B$, we write $A\ind_{C} B$ if and only if $\mathrm{acl}_{\sigma}(A)$ is algebraically independent from $\mathrm{acl}_{\sigma}(B)$ over $\mathrm{acl}_{\sigma}(C)$.
\end{definition}
\begin{fact}(See \cite[(1.9)]{acfa}.)
    The $\mathcal{L}_{\sigma}$-theory $\mathrm{ACFA}$ is a simple theory and forking independence coincides with the above notion of independence. In particular 3-amalgamation holds over $\mathrm{acl}_{\sigma}$-closed sets. 
\end{fact}

In \cite[Section 6.3]{Hrushovski2021AxsTW} Hrushovski observes that in $\mathrm{ACFA}^{+}$ 3-amalgamation may fail over $\mathrm{acl}_{\sigma}$-closed sets and conjectures that a complete characterisation of this phenomena can be obtained. This characterisation is obtained in \cite[Theorem 4.18]{ludwig2025modeltheorydifferencefields}. To state the result we recall the notion of an $\sigma$-AS-closed set.

\begin{definition}\label{torsordefinition}(See \cite[Definition 4.10]{ludwig2025modeltheorydifferencefields}.)
Given $a\in K$ we denote by $\mathfrak{T}_{a}$ the additive $F$-Torsor defined by $\sigma(x)-x=a$.
\end{definition}
\begin{definition}\label{definitiontorsorclosed}(See \cite[Definition 4.11]{ludwig2025modeltheorydifferencefields}.)
We define a set $E$ as \textit{$\sigma$-AS-closed} (for Artin-Schreier) if $E=\mathrm{acl}_{\sigma}(E)$ and, moreover, if for every $a\in E$ there exists $b\in E$ such that $b\in \mathfrak{T}_{a}$.
\end{definition}
We sketch why $E$ being $\sigma$-AS-closed is indeed a necessary condition for 3-amalgamation to hold over $E$. The idea is essentially to consider independent solutions $\alpha_{1},\alpha_{2},\alpha_{3}$ of $\mathfrak{T}_{a}$ where $a\in E$ is such that $\mathfrak{T}_{a}\cap E=\emptyset$. In this case one shows using Fact \ref{factlinindconsistent} that for any $1\leq i<j\leq3$ and $r_{ij}\in S^{1}$, it is consistent with $\mathrm{tp}(\alpha_{i}/E)\cup \mathrm{tp}(\alpha_{j}/E)$ that $\Psi(\alpha_{i}-\alpha_{j})=r_{ij}$. Using this one constructs the types extending the $\mathrm{tp}(\alpha_{i}/E)$ to a 3-amalgamation problem that does not have a solution. So, essentially, 3-amalgamation does not hold due to the presence of the following hyperimaginaries $\mathfrak{T}_{a}/E_{a}$ for $a\in E$ with $\mathfrak{T}_{a}\cap E=\emptyset$.

\begin{definition}\label{definitionequivrelationEa}(See \cite[Definition 4.14]{ludwig2025modeltheorydifferencefields}.)
    Fix $K\models \mathrm{ACFA}^{+}$ and $a\in K$. Let $E_{a}$ be the equivalence relation on the $(F,+)$-Torsor $\mathfrak{T}_{a}$ given by
    \[x\,E_{a}\,y\;\iff\;\Psi(x-y)=1.\]
    We call any $\mathbf{b}_{a}\in \mathfrak{T}_{a}/E_{a}$  a $\sigma$-AS-\textit{imaginary}.
\end{definition}
Note that the above defined hyperimaginary $\mathfrak{T}_{a}/E_{a}$ is an imaginary in the sense of continuous logic (CL-imaginary) as in \cite[Section 6]{mtfms}.

\begin{fact}(See \cite[Theorem 4.18]{ludwig2025modeltheorydifferencefields}.)\label{theorem3amalgamationcharacterisation}
    Over $E=\mathrm{acl}_{\sigma}(E)$ 3-amalgamation holds if and only if $E$ is $\sigma$-AS-closed.
\end{fact}

The statement of the theorem is then that the above described obstruction to 3-amalgamation turns out to be the only obstruction in $\mathrm{ACFA}^{+}$. The proof can be done in two main steps. First, one shows that whenever for two independent $\mathrm{acl}_{\sigma}$-closed difference fields extensions $A_{1}$ and $A_{2}$ of $E$, the following equality of $\mathbb{Q}$-vector spaces $\langle A_{1},A_{2}\rangle_{\mathbb{Q}}\cap F=\langle (A_{1}\cap F),(A_{2}\cap F)\rangle_{\mathbb{Q}}$ holds, then 3-amalgamation holds over $E$. Next, one employs the following lemma. 
\begin{lemma}\label{lemmatorsorclosedgivesvectorspace}
    If $E=\mathrm{acl}_{\sigma}(E)$ is $\sigma$-AS-closed and $A_{1}$, $A_{2}$ are independent difference field extensions of $E$, then
    \[\langle A_{1},A_{2}\rangle_{\mathbb{Q}}\cap F=\langle (A_{1}\cap F),(A_{2}\cap F)\rangle_{\mathbb{Q}}\]
\end{lemma}
\begin{proof}
 Let $a_{1}\in A_{1}$, $a_{2}\in A_{2}$ and assume $a_{1}-a_{2}\in F$. It follows by independence that $\sigma(a_{1})-a_{1}=\sigma(a_{2})-a_{2}=e\in E$. Take $\Tilde{e}\in \mathfrak{T}_{e}\cap E$, then $a_{1}-a_{2}=(a_{1}-\Tilde{e})-(a_{2}-\Tilde{e})$ and $a_{1}-\Tilde{e},a_{2}-\Tilde{e}\in F$.
\end{proof}

Next, we state the following Fact which not only allows to construct sets $E=\mathrm{acl}_{\sigma}(E)$ that are not $\sigma$-AS-closed but will also allow us to employ stability-theoretic arguments in the following sections.
Note that the proof is completely elementary. It is implied by the more general Lemma \ref{lemmalineardiffequationnotrealisedinalgclosure} which we will prove in Section \ref{sectioncounterexample}.

\begin{lemma}\label{lemmatorsornotrealisedinalgclosure}(See \cite[Lemma 4.12]{ludwig2025modeltheorydifferencefields})
Let $(K,\sigma)$ be a difference field of characteristic $0$. Assume that $\mathfrak{T}_{a}$ is not realised in $K$ where $a\in K$, then $\mathfrak{T}_{a}$ is also not realised in $K^{\mathrm{alg}}$.    
\end{lemma}

Finally let us recall that as an immediate consequence of Fact \ref{theorem3amalgamationcharacterisation} $\mathrm{ACFA}^{+}$ is simple.

\begin{fact}
    The theory $\mathrm{ACFA}^{+}$ is simple and forking independence coincides with the notion of independence given in Definition \ref{definitionindependenceacfa}.
\end{fact}

\section{A criterion for higher amalgamation in $\mathrm{ACFA}^{+}$}\label{sectioncriterionhigheramalg}

\begin{convention}In the following we will always work in some monster model $\mathcal{M}\models \mathrm{ACFA}$, if not specified otherwise. In particular independent n-systems $(\bar{a}_{w})_{w\in\mathcal{P}(\mathbf{n})}$ and the corresponding base set $E$ will be living in $\mathcal{M}$. As the independence relation does not change when passing from a model of $\mathrm{ACFA}^{+}$ to its $\mathcal{L}_{\sigma}$-reduct, our results (when they only consider the $\mathcal{L}_{\sigma}$-structure) will be applicable when working in an ambient model of $\mathrm{ACFA}^{+}$ as well. Hence, we will sometimes not explicitly differentiate between working with an independent n-system in $\mathrm{ACFA}$ or $\mathrm{ACFA}^{+}$.
\end{convention}

The idea of this section will be to develop a higher dimensional analogue (Definition \ref{definitionffequationandffdecomposable}) of the property in Lemma \ref{lemmatorsorclosedgivesvectorspace}and then to show that it holds (over $E$) if and only if $n$-amalgamation holds (over $E$). In the subsequent sections we will then see  that to determine for which $E$ this property will hold and what the connection to the $\sigma$-AS map $\wp_{\sigma}(x):=\sigma(x)-x$ is, will, however, turn out to be surprisingly more intricate.\\
Now, to formulate the higher dimensional analogue of the property in Lemma \ref{lemmatorsorclosedgivesvectorspace}, we will, as already in \cite{ludwig2025modeltheorydifferencefields}, work on the level of additive formulas directly, as this will be useful in the proofs. In Lemma \ref{lemmacorrespondingvectorpsaceproperty} we then state the corresponding property on the level of vector spaces.

\begin{definition}\label{definitionffequationandffdecomposable}
Let $(\bar{a}_{w})_{w\in\mathcal{P}(\mathbf{n})}$ be an independent n-system over $E$.   We call an equation of the form $\sum_{1\leq i\leq n}b_{\widehat{w}_{i}}=0$ a \textit{fixed field additive equation of height $n$ over $E$} (or simply \textit{ff-additive equation}) if for every $1\leq i\leq n$ we have $b_{\widehat{w}_{i}}\in F\cap \bar{a}_{\widehat{w}_{i}}$. If we want to specify the ambient independent $n$-system over $E$, we say that the equation is \textit{in} $(\bar{a}_{w})_{w\in\mathcal{P}(\mathbf{n})}$.
We call the equation \textit{ff-decomposable} if it is decomposable into fixed field elements, i.e., (working with the notation from Definition \ref{definitiondecomposableequation}) if we find a decomposition with $c_{\widehat{w}_{i}}^{j}\in F\cap \bar{a}_{\widehat{w}_{i,j}}$ for any pair $1\leq i,j\leq n$ with $i\neq j$.
\end{definition}

To prove Lemma \ref{lemmanotffdecomposableyieldsnotamalgamation} we will make use of the following lemma that recovers step 2 of the proof of Lemma \ref{lemmadecomposabilityofgeneraladditiveequations} in the context of an ff-additive equation.
\begin{lemma}\label{lemmarecoveringstep2ff}
Fix an independent $n$-system $(\bar{a}_{w})_{w\in\mathcal{P}(\mathbf{n})}$ over some $E=\mathrm{acl}_{\sigma}(E)$.
Assume that either $n=3$ or $n>3$ and for any $1\leq i\leq n$ every ff-additive equation in $(\Tilde{a}_{w})_{w\in\mathcal{P}(\widehat{\mathbf{n}}_{i})}$ over $\bar{a}_{i}$ (using Notation \ref{remarkindependentsystemonesteplower}) of height $n-1$ is ff-decomposable over $\bar{a}_{i}$.
Let $\sum_{1\leq i\leq n}b_{\widehat{w}_{i}}=0$ be an ff-additive equation over $E$ of height $n$ and assume that for any $1\leq i\leq n$ we have that $b_{\widehat{w}_{i}}=\sum_{j\neq i}d_{\widehat{w}_{i}}^{j}$ where $d_{\widehat{w}_{i}}^{j}\in F\cap \bar{a}_{\widehat{w}_{i,j}}$ (but not necessarily $d_{\widehat{w}_{i}}^{j}=-d_{\widehat{w}_{j}}^{i}$). Then $\sum_{1\leq i\leq n}b_{\widehat{w}_{i}}=0$ is even ff-decomposable.
\end{lemma}
\begin{proof}
This follows using the exact same argument as in step 2 of the proof of Lemma \ref{lemmadecomposabilityofgeneraladditiveequations}. Concretely, if $n=3$ we note that $\delta_{i}=d_{\widehat{w}_{j}}^{k}-d_{\widehat{w}_{k}}^{j}$ for $\{i,k,j\}=\{1,2,3\}$ are elements of the fixed field, since the $d_{\widehat{w}_{j}}^{k}$ are. Hence, the same procedure as in step 2 of the proof of Lemma \ref{lemmadecomposabilityofgeneraladditiveequations} yields a system of $c_{\widehat{w}_{i}}^{j}\in F$ witnessing ff-decomposability.\\
If $n>3$, we proceed as in the induction step of the proof of Lemma \ref{lemmadecomposabilityofgeneraladditiveequations}: We apply the assumption that $b_{\widehat{w}_{n}}=\sum_{j\neq n}d_{\widehat{w}_{n}}^{j}$ where $d_{\widehat{w}_{n}}^{j}\in F\cap \bar{a}_{\widehat{w}_{n,j}}$ and by the exact same argument as in the induction step of the proof of Lemma \ref{lemmadecomposabilityofgeneraladditiveequations} we obtain an ff-additive equation of height $n-1$ over $\bar{a}_{n}$, which is then ff-decomposable. In this way again, we obtain a system of $c_{\widehat{w}_{i}}^{j}$ that is given by fixed field elements.
\end{proof}

\begin{lemma}\label{lemmaffdecomposableimpliesnamalgam}
    Let $n\geq 3$ and $E=\mathrm{acl}_{\sigma}(E)\subseteq \mathcal{M}\models \mathrm{ACFA}^{+}$. If for every $3\leq m\leq n$ every ff-additive equation of height $m$ is ff-decomposable over $E$, then $\mathrm{Th}(M)$ has n-amalgamation over $E$.
\end{lemma}
\begin{proof}
Let $(p_{w})_{w\in \mathcal{P}^{-}(\mathbf{n})}$ be an n-amalgamation problem in $\mathrm{ACFA}^{+}$ over $E$. Let $\Tilde{p}_{\mathbf{n}}$ be a solution to the n-amalgamation problem $(p_{w}\restriction_{\mathcal{L}_{\sigma}})_{w\in \mathcal{P}^{-}(\mathbf{n})}$ in $\mathrm{ACFA}$. We will extend $\Tilde{p}_{\mathbf{n}}$ to a complete type $p_{\mathbf{n}}$ in $\mathrm{ACFA}^{+}$ which completes the system $(p_{w})_{w\in \mathcal{P}^{-}(\mathbf{n})}$.
Let $a_{\mathbf{n}}$ be a realisation of $\Tilde{p}_{\mathbf{n}}$ and denote by $\bar{a}_{w}$ the subset of $a_{\mathbf{n}}$ corresponding to the realisation of $p_{w}\restriction_{\mathcal{L}_{\sigma}}$. Then $(\bar{a}_{w})_{w\in\mathcal{P}(\mathbf{n})}$ is an independent $n$-system over $E$. Next, we equip the set $A:=\bigcup_{w\in\mathcal{P}^{-}(\mathbf{n})}\bar{a}_{w}\cap F$ with a map $\Psi: A\rightarrow S^{1}$ such that $\bar{a}_{w}$ together with $\Psi\restriction_{\bar{a}_{w}\cap F}$ realises $p_{w}$. Note that by construction, all the restrictions $\Psi\restriction_{\bar{a}_{w}\cap F}:(\bar{a}_{w}\cap F, +)\rightarrow (S^{1},\cdot)$ yield group homomorphisms.\\
Let $V$ denote the $\mathbb{Q}$-vector space generated by $A$. Then, it suffices to show that $\Psi$ extends to a homomorphism on $V$ because once this is done we can choose $\Psi$ freely outside of $V$ by Fact \ref{factlinindconsistent}. By compactness $\Psi$ extends to a homomorphism on $V$ if and only if for every ff-additive equation $\sum_{1\leq i\leq n}b_{\widehat{w}_{i}}=0$ we have that $\prod_{1\leq i\leq n}\Psi(b_{\widehat{w}_{i}})=1$. Now, if $\sum_{1\leq i\leq n}b_{\widehat{w}_{i}}=0$ is ff-decomposable, we have $b_{\widehat{w}_{i}}=\sum_{1\leq j\leq n,\,j\neq i}c_{\widehat{w}_{i}}^{j}$ where $c_{\widehat{w}_{i}}^{j}\in \bar{a}_{i,j}\cap F$. Then, since $\Psi$ is a homomorphism on $\bar{a}_{\widehat{w}_{i}}\cap F$, we have for every $1\leq i\leq n$ that $\Psi(b_{\widehat{w}_{i}})=\prod_{1\leq j\leq n,\,j\neq i}\Psi(c_{\widehat{w}_{i}}^{j})$ and consequently we obtain
\[\prod_{1\leq i\leq n}\Psi(b_{\widehat{w}_{i}})=\prod_{1\leq i\leq n}\prod_{1\leq j\leq n,\,j\neq i}\Psi(c_{\widehat{w}_{i}}^{j}).\]
But since we have that $c_{\widehat{w}_{i}}^{j}=-c_{\widehat{w}_{j}}^{i}$ and since $\Psi$ is a homomorphism on $\bar{a}_{\widehat{w}_{i,j}}\cap F$ it follows that $\Psi(c_{\widehat{w}_{i}}^{j})=\Psi(c_{\widehat{w}_{j}}^{i})^{-1}$ and thus
\[\prod_{1\leq i\leq n}\Psi(b_{\widehat{w}_{i}})=\prod_{1\leq i\leq n}\prod_{1\leq j\leq n,\,j\neq i}\Psi(c_{\widehat{w}_{i}}^{j})=1.\]
\end{proof}

The following lemma concerns the opposite direction of Lemma \ref{lemmaffdecomposableimpliesnamalgam}. More concretely, we will show how to obtain from an ff-additive equation that is not ff-decomposable an amalgamation problem which does not have a solution.

\begin{lemma}\label{lemmanotffdecomposableyieldsnotamalgamation}
Let $E=\mathrm{acl}_{\sigma}(E)\subseteq\mathcal{M}\models \mathrm{ACFA}^{+}$. Assume that there is an independent $n$-system $(\bar{a}_{w})_{w\in\mathcal{P}(\mathbf{n})}$ over $E$ with an ff-additive equation $0=\sum_{1\leq i\leq n}b_{\widehat{w}_{i}}$ of height $n$ that is not ff-decomposable and $n$ is minimal as such. Then, there is an $n$-amalgamation problem $(p_{w})_{w\in \mathcal{P}^{-}(\mathbf{n})}$ over $E$ that does not have a solution such that moreover $(p_{w}\restriction_{\mathcal{L}_{\sigma}})_{w\in \mathcal{P}^{-}(\mathbf{n})}=(\mathrm{tp}_{\mathcal{L}_{\sigma}}(\bar{a}_{w}/E))_{w\in \mathcal{P}^{-}(\mathbf{n})}$.
\end{lemma}
\begin{proof}
The proof is by induction over $n$, starting with $n=3$ which is given by \cite[Lemma 4.17]{ludwig2025modeltheorydifferencefields}. As $E$ is an arbitrary $\mathrm{acl}_{\sigma}$-closed set, we will assume in the induction step that the statement is shown for $(n-1)$ for all $\mathrm{acl}_{\sigma}$-closed sets.
We define $p_{w}\restriction_{\mathcal{L_{\sigma}}}:=\mathrm{tp}(\bar{a}_{w}/E)\restriction_{\mathcal{L}_{\sigma}}$ and set $\Psi$ to be a group homomorphism (into $S^{1}$) on the $\mathbb{Q}$-vector space generated by $\bigcup_{1\leq i,j\leq n,\,i\neq j}F\cap \bar{a}_{\widehat{w}_{i,j}}$ and thus obtain the corresponding $\mathcal{L}_{\sigma}^{+}$-types $p_{\widehat{w}_{i,j}}$. The goal will now be to extend $\Psi$ to $\bar{a}_{\widehat{w}_{i}}$ for any $1\leq i\leq n$ consistent with the $p_{\widehat{w}_{i,j}}$ such that however $\Psi(b_{\widehat{w}_{k}})\neq\prod_{1\leq i\leq n,\,i\neq k}\Psi(b_{\widehat{w}_{i}})$ would have to hold and thus yielding a counterexample to $n$-amalgamation.\\
We will first show that the first part of the assumption of Lemma \ref{lemmarecoveringstep2ff} is satisfied, i.e., we show that for any $1\leq i\leq n$, every ff-additive equation in $(\Tilde{a}_{w})_{w\in\mathcal{P}(\widehat{\mathbf{n}}_{i})}$ over $\bar{a}_{i}$ of height $n-1$ is ff-decomposable over $\bar{a}_{i}$.
Assume otherwise, then for some $1\leq i\leq n$,
we can use the induction hypothesis (working over the $\mathrm{acl}_{\sigma}$-closed set $\bar{a}_{\widehat{w}_{i}}$), to obtain $(\Tilde{p}_{w})_{w\in\mathcal{P}^{-}(\widehat{\mathbf{n}}_{i})}$, a counterexample to $(n-1)$-amalgamation over $\bar{a}_{i}$ which extends $(\mathrm{tp}_{\mathcal{L}_{\sigma}}(\Tilde{a}_{w})/\bar{a}_{i})_{w\in\mathcal{P}^{-}(\widehat{\mathbf{n}}_{i})}$.\\
Note that $(\Tilde{p}_{w})_{w\in\mathcal{P}^{-}(\widehat{\mathbf{n}}_{i})}$ already induces a partial $n$-amalgamation problem over $E$ (with the $\widehat{\mathbf{n}}_{i}$-component missing) as $\bar{a}_{i}$ is independent from $\bar{a}_{\widehat{w}_{i}}$ over $E$. This even extends to a full $n$-amalgamation problem over $E$ by amalgamating the types corresponding to $w\in \widehat{\mathbf{n}}_{i}$ (as otherwise $E$ would not have $(n-1)$-amalgamation which contradicts the assumption by Lemma \ref{lemmaffdecomposableimpliesnamalgam}).\\
Hence, since the first assumption of Lemma \ref{lemmarecoveringstep2ff} is satisfied and moreover $\sum_{1\leq i\leq n}b_{\widehat{w}_{i}}=0$ is not ff-decomposable, it follows from Lemma \ref{lemmarecoveringstep2ff} that there is some $1\leq k\leq n$ such that $b_{\widehat{w}_{k}}$ can not be written as $b_{\widehat{w}_{k}}=\sum_{j\neq k}d_{\widehat{w}_{k}}^{j}$ for $d_{\widehat{w}_{k}}^{j}\in F\cap \bar{a}_{\widehat{w}_{k,j}}$. 
This means that $b_{\widehat{w}_{k}}$ is not contained in the $\mathbb{Q}$-vector space generated by all the $\bar{a}_{\widehat{w}_{i,j}}\cap F$ where $1\leq j\leq n,\,k\neq j$ and thus for any $r\in S^{1}$ it follows that $r=\Psi(b_{\widehat{w}_{k}})$ is consistent with $\bigcup_{1\leq j\leq n,\,k\neq j}p_{\widehat{w}_{i,k}}$ by Fact \ref{factlinindconsistent}.
So, we can extend $\Psi$ to a homomorphism on $F\cap \bar{a}_{\widehat{w}_{i}}$ for $1\leq i\leq n$ such that $\Psi(b_{\widehat{w}_{k}})\neq\prod_{1\leq i\leq n,\,i\neq k}\Psi(b_{\widehat{w}_{i}})$ would have to hold. For every $1\leq i\leq n$ we set $p_{\widehat{w}_{i}}:=\mathrm{tp}(\bar{a}_{\widehat{w}_{i}})$. Those types are consistent with the $p_{\widehat{w}_{i,j}}$
and hence yield an $n$-amalgamation over $E$ that does not have a solution.
\end{proof}

As mentioned before we can give an equivalent version for vector spaces of our criterion for $n$-amalgamation. We will not use it later but give it for the purpose of illustration.
\begin{lemma}\label{lemmacorrespondingvectorpsaceproperty}
 Every ff-additive equation of height $3\leq k\leq n$ over $E=\mathrm{acl}_{\sigma}(E)$ is ff-decomposable if and only if for all $2\leq m\leq n-1$ and every independent m-system $(\bar{a}_{w})_{w\in\mathcal{P}(\mathbf{m})}$ (in $\mathcal{L}_{\sigma}$) over $E$ the following property holds
 \[\langle\bar{a}_{\widehat{w}_{1}},\dots,\bar{a}_{\widehat{w}_{m}}\rangle_{\mathbb{Q}}\cap F=\langle  \bar{a}_{\widehat{w}_{1}}\cap F,\dots,\bar{a}_{\widehat{w}_{m}}\cap F \rangle_{\mathbb{Q}}.\]
\end{lemma}
\begin{proof}
One direction is trivial. The other one, that the condition on vector spaces yields the one about ff-decomposability, follows from Lemma \ref{lemmarecoveringstep2ff}.
\end{proof}

\begin{remark}
In \cite{grupoidshrushovski} Hrushovski establishes for a stable theory $T$ a connection between $n$-amalgamation and $n-1$-uniqueness (for $n\geq 3$). Here, $n$-uniqueness will mean that any $n$-amalgamation problem has a unique solution. More precisely, Hrushovski shows that for $T$ stable, $4$-amalgamation is equivalent to 3-uniqueness. It was then proved in \cite{typeamalgamationandpolygroupoids} that the same argument works for $n$-amalgamation and $n-1$-uniqueness, if we not only assume stability but also $m$-uniqueness for all $2\leq m\leq n-2$.\\
In our context, however, $n$-uniqueness already fails in the case $n=2$ ($\mathrm{ACFA}$ is not stable). So, the above correspondence cannot be applied directly. Let us note in passing that Lemma \ref{lemmacorrespondingvectorpsaceproperty} describes \textit{some} version of $n-1$-uniqueness, namely the uniqueness of extensions of $\Psi$ to the $\mathbb{Q}$-vector space given by $F$ in a solution to an $n-1$-amalgamation problem.
\end{remark}

\section{Higher amalgamation over models}\label{sectionnamalgovermodelschapterfour}
In this section we will prove that $n$-amalgamation holds for all $n\in\mathbb{N}$ over substructures that are models of $\mathrm{ACFA}$, so in particular over models of $\mathrm{ACFA}^{+}$. Unlike in $\mathrm{ACFA}$ or $\mathrm{PF}^{+}$ where $n$-amalgamation holds over all algebraically closed sets, in $\mathrm{ACFA}^{+}$ this is not a straight-forward generalisation of the results on 3-amalgamation. On the contrary, if $\mathrm{ACFA}^{+}$ was a rank 1 theory (which it is not), it would even be easy to construct a counterexample to 4-amalgamation over a model using the non-trivial obstruction to 3-amalgamation. More concretely, we would simply take some $a$ and model $\mathcal{N}\models \mathrm{ACFA}^{+}$ such that $\mathfrak{T}_{a}$ is not realised in $\mathrm{acl}_{\sigma}(Na)$. Then, we construct a 3-amalgamation problem induced by $\{\bar{a}_{1},\bar{a}_{2},\bar{a}_{3}\}$ over $\bar{a}_{4}:=\mathrm{acl}_{\sigma}(Na)$ which does not have a solution over $\bar{a}_{4}$ with $a_{i}\in\mathfrak{T}_{a}\cap\bar{a}_{i}$ for $1\leq i\leq 3$. Next, we consider the system of types induced by $\{\bar{a}_{1},\bar{a}_{2},\bar{a}_{3},\bar{a}_{4}\}$ and if this was a 4-amalgamation problem over $N$ we would indeed have found a counterexample to 4-amalgamation over a model. Now, however unlike in the rank 1 case the above system does not yield a 4-amalgamation problem in $\mathrm{ACFA}^{+}$: Recall that the 3-amalgamation problem was constructed using realisations $a_{i}\in\bar{a}_{i}$ of $\mathfrak{T}_{a}$, so in particular $\sigma(a_{i})-a_{i}=a$ prohibits independence of $\bar{a}_{i}$ and $\bar{a}_{4}$ for $1\leq i\leq 3$. We will now see that indeed n-amalgamation holds over models. While some of the ideas from the context of 3-amalgamation are present in the proof, it is not simply a straight-forward higher dimensional adaptation and will take up the rest of this section. We start with a definition that can be seen as a higher dimensional analogue of being $\sigma$-AS-closed.
\begin{definition}\label{definitionntorsorclosed}
    We say that some set $E\subseteq K\models \mathrm{ACFA}$ is \textit{n-$\sigma$-AS-closed} if $E=\mathrm{acl}_{\sigma}(E)$ and the following holds for all $1\leq m\leq n$ and any independent $m$-system $(\bar{a}_{w})_{w\in\mathcal{P}(\mathbf{m})}$ over $E$: For any $\Tilde{b}=\sum_{1\leq i\leq m}\Tilde{b}_{\widehat{w}_{i}}$ with $\Tilde{b}_{\widehat{w}_{i}}\in \bar{a}_{\widehat{w}_{i}}$ it follows that whenever $\mathfrak{T}_{\Tilde{b}}$ is realised in $a_{\mathbf{m}}$, then there exist some $b_{\widehat{w}_{i}}\in \bar{a}_{\widehat{w}_{i}}$ such that $\Tilde{b}=\sum_{1\leq i\leq m}b_{\widehat{w}_{i}}$ and
    $\mathfrak{T}_{b_{\widehat{w}_{i}}}$ is realised in $\bar{a}_{\widehat{w}_{i}}$ for all $1\leq i\leq m$.
\end{definition}
\begin{notation}\label{notationoldversionwitnesstorsorclosed}
    If for some specific independent $n$-system $(\bar{a}_{w})_{w\in\mathcal{P}(\mathbf{n})}$ the above holds for all subsystems $(\bar{a}_{w})_{w\in\mathcal{P}(\mathbf{n}\backslash J)}$ where $J\subsetneq\mathbf{n}$, then we say that $E$ is $n$-$\sigma$-AS-closed for $(\bar{a}_{w})_{w\in\mathcal{P}(\mathbf{n})}$.
    On the contrary, if there is some $\Tilde{b}$ as above such that $\mathfrak{T}_{\Tilde{b}}$ is realised in $a_{\mathbf{n}}$ but for any sum $\Tilde{b}=\sum_{1\leq i\leq n}b_{\widehat{w}_{i}}$ some $\mathfrak{T}_{b_{\widehat{w}_{i}}}$ is not realised in $\bar{a}_{\widehat{w}_{i}}$, we say that $(\bar{a}_{w})_{w\in\mathcal{P}(\mathbf{n})}$ witnesses that $E$ is not $n$-$\sigma$-AS-closed.
\end{notation}

\begin{remark}
 The set $E$ is $1$-$\sigma$-AS-closed if and only if $E$ is $\sigma$-AS-closed in the sense of Definition \ref{definitiontorsorclosed}.
\end{remark}

\begin{lemma}\label{lemmantorsorclosednessonestepbelow}
    Let $n\geq 2$ and $E=\mathrm{acl}_{\sigma}(E)$ be $n$-$\sigma$-AS-closed for $(\bar{a}_{w})_{w\in\mathcal{P}(\mathbf{n})}$, an independent $n$-system over $E$. As in Notation \ref{remarkindependentsystemonesteplower} we define for every $w\in \mathcal{P}(\mathbf{n-1})$ the tuple $\Tilde{a}_{w}:=\mathrm{acl}_{\sigma}(\bar{a}_{w}\cup \bar{a}_{n})=\bar{a}_{w\cup\{n\}}$.
    Then, $\bar{a}_{n}\supset E$ is $(n-1)$-$\sigma$-AS-closed for $(\Tilde{a}_{w})_{w\in\mathcal{P}(\mathbf{n-1})}$.
\end{lemma}
\begin{proof}
We consider a sum $\Tilde{b}=\sum_{1\leq i\leq n-1}\Tilde{b}_{\widehat{w}_{i}}$ in the independent $(n-1)$-system $(\Tilde{a}_{w})_{w\in\mathcal{P}(\mathbf{n-1})}$ over $\bar{a}_{n}$. In particular $\Tilde{b}_{\widehat{w}_{i}}\in \bar{a}_{\widehat{w}_{i}\cup\{n\}}$ holds for every $1\leq i\leq n-1$. We now work over the set $\mathbf{n}$ and define $\Tilde{b}_{\widehat{w}_{n}}:=0$. We consider the sum  $\Tilde{b}=\sum_{1\leq i\leq n}\Tilde{b}_{\widehat{w}_{i}}$ in $(\bar{a}_{w})_{w\in\mathcal{P}(\mathbf{n})}$. By $n$-$\sigma$-AS-closedness we obtain some $b_{\widehat{w}_{i}}\in \bar{a}_{\widehat{w}_{i}}$ such that \[\Tilde{b}=\sum_{1\leq i\leq n-1}\Tilde{b}_{\widehat{w}_{i}}=\sum_{1\leq i\leq n}b_{\widehat{w}_{i}}\;\;\;(\star)\] and, moreover, $\mathfrak{T}_{b_{\widehat{w}_{i}}}$ is realised in $\bar{a}_{\widehat{w}_{i}}$ for any $1\leq i\leq n$. If $b_{\widehat{w}_{n}}=0$ holds, then the statement already follows. Hence, we assume $b_{\widehat{w}_{n}}\neq0$.\\
From $(\star)$ we obtain $b_{\widehat{w}_{n}}=\sum_{1\leq i\leq n-1}(\Tilde{b}_{\widehat{w}_{i}}-b_{\widehat{w}_{i}})$. We can apply Lemma \ref{lemmadecomposabilityofgeneraladditiveequations} to find $c_{\widehat{w}_{i,n}}\in \bar{a}_{\widehat{w}_{i,n}}$ such that $b_{\widehat{w}_{n}}=\sum_{1\leq i\leq n-1}c_{\widehat{w}_{i,n}}.$ As $\mathfrak{T}_{b_{\widehat{w}_{n}}}$ is realised in $\bar{a}_{\widehat{w}_{n}}$ (and by assumption $E$ is $(n-1)$-$\sigma$-AS-closed for $(\bar{a}_{w})_{w\in\mathcal{P}(\widehat{\mathbf{n}}_{n})}$) we find $d_{\widehat{w}_{i,n}}\in\bar{a}_{\widehat{w}_{i,n}}$ such that $b_{\widehat{w}_{n}}=\sum_{1\leq i\leq n-1}d_{\widehat{w}_{i,n}}$ and $\mathfrak{T}_{d_{\widehat{w}_{i,n}}}$ is realised in $\bar{a}_{\widehat{w}_{i,n}}$ for all $1\leq i\leq n-1$. We insert this in $(\star)$ and obtain
\[\Tilde{b}=\sum_{1\leq i\leq n-1}\Tilde{b}_{\widehat{w}_{i}}=\sum_{1\leq i\leq n}b_{\widehat{w}_{i}}=\sum_{1\leq i\leq n-1}(b_{\widehat{w}_{i}}+d_{\widehat{w}_{i,n}})\]
For all $1\leq i\leq n-1$ we set $e_{\widehat{w}_{i}}:=b_{\widehat{w}_{i}}+d_{\widehat{w}_{i}}$. Since for any $1\leq i\leq n-1$ both $\mathfrak{T}_{b_{\widehat{w}_{i}}}$ and $\mathfrak{T}_{d_{\widehat{w}_{i,n}}}$ are realised in $\bar{a}_{\widehat{w}_{i}}$, it follows that $\mathfrak{T}_{e_{\widehat{w}_{i}}}$ is realised in $\bar{a}_{\widehat{w}_{i}}=\Tilde{a}_{\widehat{w}_{i}}$ where in the last expression we consider $\widehat{w}_{i}$ as a subset of $\mathbf{n-1}$. This completes the proof.
\end{proof}
The next goal will be to generalise Theorem \ref{theorem3amalgamationcharacterisation} to $n>3$. Before we do so, we prove that a model of $\mathrm{ACFA}$ is indeed $n$-$\sigma$-AS-closed for any $n\in\mathbb{N}$. From this we will then be able to deduce that $n$-amalgamation in $\mathrm{ACFA}^{+}$ holds over substructures that are models of $\mathrm{ACFA}$. The next lemma heavily relies on Lemma \ref{lemmatorsornotrealisedinalgclosure}. Once Lemma \ref{lemmatorsornotrealisedinalgclosure} is applied the argument only uses stability of $\mathrm{ACF}$ and resembles a lot the proof of Lemma \ref{lemmadecomposabilityofgeneraladditiveequations}. 

\begin{lemma}\label{lemmaspecialisationtorsorformula}
    Let $\gamma(x_{0},\dots,x_{k})$ be an $\mathcal{L}_{\mathrm{ring}}$-term. The formula $\phi(x_{0},\dots,x_{k})$, given by
    \[\exists z\;\sigma(z)-z=\gamma(x_{0},\dots,x_{k}),\]
    specialises over any $E\models \mathrm{ACFA}$.
\end{lemma}
\begin{proof}
We fix an independent $n$-system $(\bar{a}_{w})_{w\in\mathcal{P}(\mathbf{n})}$ over $E$ and we assume that $\phi(b_{J_{0}},\dots,b_{J_{k}})$ holds (in $\bar{a}_{\mathbf{n}}$) for elements $b_{J_{i}}\in \bar{a}_{J_{i}}$ and $J_{0},\dots,J_{k}\subset \{1,\dots,n\}$. We follow the structure and notation of the proof of Lemma \ref{lemmadecomposabilityofgeneraladditiveequations}. The goal will be to specialise $\phi(b_{J_{0}},\dots,b_{J_{k}})$ to $\bar{a}_{J_{0}}$. By Lemma \ref{lemmatorsornotrealisedinalgclosure} we know that there is some $c$ in the field composite of $\bar{a}_{J_{0}}$ and $\bar{a}_{\widehat{w}_{J_{0}}}$ (as by assumption there is one in $\bar{a}_{\mathbf{n}}$, the algebraic closure of $\bar{a}_{J_{0}}$ and $\bar{a}_{\widehat{w}_{J_{0}}}$) such that $\sigma(c)-c=\gamma(b_{J_{0}},\dots,b_{J_{k}})$ holds. Hence, there are finite tuples $\bar{d}=d_{1},\dots,d_{n}$ in $\bar{a}_{J_{0}}$ and $\bar{e}=e_{1},\dots,e_{m}$ in $\bar{a}_{\widehat{w}_{J_{0}}}$ as well as $\mathcal{L}_{\mathrm{ring}}$-terms $f(\bar{y},\bar{z})$and $g(\bar{y},\bar{z})$ such that $c=\frac{f(\bar{d},\bar{e})}{g(\bar{d},\bar{e})}$ and consequently $\sigma(c)=\frac{f(\sigma(\bar{d}),\sigma(\bar{e}))}{g(\sigma(\bar{d}),\sigma(\bar{e}))}$.\\
   As in the proof of Lemma \ref{lemmadecomposabilityofgeneraladditiveequations} we find for any $1\leq i\leq k$ an $\mathcal{L}_{\mathrm{ring}}(\bar{a}_{J_{i}\cap J_{0}})$-formula $q_{i}(x_{i},\bar{w}_{i})$ such that there is some tuple $\bar{\beta}_{i}\in \bar{a}_{J_{i}\backslash J_{0}}$ such that $q_{i}(b_{J_{i}},\bar{\beta}_{i})$ expresses that $b_{J_{i}}$ is in the algebraic closure of the field composite of $\bar{a}_{J_{i}\cap J_{0}}$ and $\bar{a}_{J_{i}\backslash J_{0}}$. We set $\bar{\beta}=(\bar{\beta}_{1},\dots,\bar{\beta}_{k})$ which by definition is a tuple in $\bar{a}_{\widehat{w}_{J_{0}}}$. Next, we define the $\mathcal{L}_{\mathrm{ring}}(\bar{a}_{J_{0}})$-formula $\rho(\bar{x},\bar{w})$ by
  \[\rho(\bar{x},\bar{w})=\bigwedge_{1\leq i\leq k}q_{i}(x_{i},\bar{w}_{i}).\]
  Consider the $\mathcal{L}_{\mathrm{ring}}(\bar{a}_{J_{0}})$-formula $\psi(\bar{u})=\psi(\bar{z},\bar{v},\bar{w},x_{1},\dots,x_{k})$ which expresses that
    \[\rho(\bar{x},\bar{w})\;\;\land\;\;\frac{f(\sigma(\bar{d}),\bar{z})}{g(\sigma(\bar{d}),\bar{z})}-\frac{f(\bar{d},\bar{v})}{g(\bar{d},\bar{v})}=\gamma(b_{J_{0}},x_{1},\dots,x_{k}).\]
Let $\mathrm{tp}_{0}(\bar{a}_{J_{0}}/E)$ denote $\mathrm{tp}(\bar{a}_{J_{0}}/E)\restriction_{\mathcal{L}_{\mathrm{ring}}}$.  By stability of $\mathrm{ACF}$, i.e., definability of the type $\mathrm{tp}_{0}(\bar{a}_{J_{0}}/E)$, there is an $\mathcal{L}_{\mathrm{ring}}(E)$-formula $\Tilde{\psi}(\bar{u})$
such that $\psi(E^{|\bar{u}|})=\Tilde{\psi}(E^{|\bar{u}|})$ holds. We note that $\mathrm{tp}_{0}(\bar{a}_{J_{0}}/\bar{a}_{\widehat{w}_{J_{0}}})$ is the unique non-forking extension of $\mathrm{tp}_{0}(\bar{a}_{J_{0}}/E)$ and consequently $\mathcal{M}\models \Tilde{\psi}(\bar{e},\sigma(\bar{e}),\bar{\beta},b_{J_{1}},\dots,b_{J_{k}})$. Hence, \[\mathcal{M}\models \exists \bar{z}\bar{v}\bar{w}\bar{x}\;\;\left(\Tilde{\psi}(\bar{z},\bar{v},\bar{w},\bar{x})\;\land\;\sigma(\bar{z})=\bar{v}\right).\] Since $E$ is an existentially closed submodel of $\mathcal{M}$ the same holds in $E$ which completes the proof.
\end{proof}

\begin{proposition}\label{propositionmodelsarentorsorclosed}
    Any $E\models \mathrm{ACFA}$ is $n$-$\sigma$-AS-closed for every $n\in\mathbb{N}$.
\end{proposition}
\begin{proof}
The proof is by induction. However, in order to make proper use of the induction hypothesis, we will have to prove something slightly stronger. We show $n$-$\sigma$-AS-closedness for an independent $m$-system $(\bar{a}_{w})_{w\in\mathcal{P}(\mathbf{m})}$ over $E$ where $m\geq n$ is allowed. (This will be made precise below.) This is in direct analogy with the notion of $(m,n)$-amalgamation from Definition \ref{definitionn-m-amalgam}.
Precisely, we will show the following.\\
\textit{Let $H:=\{J\subseteq\{1,\dots,m\}\;|\;m-|J|=n-1\}$. If for some given $\Tilde{b}=\sum_{J\in H}\Tilde{b}_{\widehat{w}_{J}}$ with $\Tilde{b}_{\widehat{w}_{J}}\in\bar{a}_{\widehat{w}_{J}}$ there is a realisation of $\mathfrak{T}_{\Tilde{b}}$ in $\bar{a}_{\mathbf{m}}$, then there are $b_{\widehat{w}_{J}}\in \bar{a}_{\widehat{w}_{J}}$, such that $\Tilde{b}=\sum_{J\in H}b_{\widehat{w}_{J}}$ and $\mathfrak{T}_{b_{\widehat{w}_{J}}}$ is realised in $\bar{a}_{\widehat{w}_{J}}$ for any $J\in H$.}\\
    First note that the case $n=1$ is still covered by the fact that $K$ is $\sigma$-AS-closed. For the induction step, we will use the following two claims.
    \begin{itemize}
        \item \textbf{Step 1}:\textit{ For any $J\in H$ the following holds:  For any $k\notin  J$ we find $c_{\widehat{w}_{J}}^{k}\in\bar{a}_{\widehat{w}_{J\cup\{k\}}}$ such that for $\Tilde{c}_{\widehat{w}_{J}}:=\sum_{k\notin J}c_{\widehat{w}_{J}}^{k}$, the set $\mathfrak{T}_{\Tilde{b}_{\widehat{w}_{J}}+\Tilde{c}_{\widehat{w}_{J}}}$ is realised in $\bar{a}_{\widehat{w}_{J}}$.}
        \item \textbf{Step 2}:\textit{ Let the $\Tilde{c}_{\widehat{w}_{J}}$ be the same as in Step 1. For $\Tilde{c}=\sum_{J\in H}\Tilde{c}_{\widehat{w}_{J}}$ there are $g_{\widehat{w}_{J}}^{k}\in\bar{a}_{\widehat{w}_{J\cup\{k\}}}$ such that \[\Tilde{c}=\sum_{J\in H,k\notin J}g_{\widehat{w}_{J}}^{k}\] and any $\mathfrak{T}_{g_{\widehat{w}_{J}}^{k}}$ is realised in $\bar{a}_{\widehat{w}_{J\cup\{k\}}}$.}
    \end{itemize}
Before proving the two claims, let us show how to deduce the statement of the lemma from them: Define
\[b_{\widehat{w}_{J}}:=\Tilde{b}_{\widehat{w}_{J}}+\Tilde{c}_{\widehat{w}_{J}}-\sum_{k\notin J}g_{\widehat{w}_{J}}^{k}.\] 
    It follows that $\sum_{J\in H}b_{\widehat{w}_{J}}=\left(\sum_{J\in H}\Tilde{b}_{\widehat{w}_{J}}\right)+\Tilde{c}-\Tilde{c}=\Tilde{b}$. Moreover, $\mathfrak{T}_{\Tilde{b}_{\widehat{w}_{J}}+\Tilde{c}_{\widehat{w}_{J}}}$ is realised in $\bar{a}_{\widehat{w}_{J}}$ for any $J\in H$. Also, for any $J,k\in H\backslash J$ there is a realisation of $\mathfrak{T}_{-g_{\widehat{w}_{J}}^{k}}$ in $\bar{a}_{\widehat{w}_{J\cup\{k\}}}\subset \bar{a}_{\widehat{w}_{J}}$. Hence, it follows that $\mathfrak{T}_{b_{\widehat{w}_{J}}}$ is realised in $\bar{a}_{\widehat{w}_{J}}$ which completes the proof.\\
    \textit{Proof of Claim 1.}  We consider the formula \[\phi(\bar{x})\equiv \exists z\;\sigma(z)-z=\sum_{J\in H}x_{\widehat{w}_{J}}.\]
    Then $\bar{a}_{\mathbf{m}}\models\phi\left((\Tilde{b}_{\widehat{w}_{J}})_{J\in H}\right)$ holds and we can apply Lemma \ref{lemmaspecialisationtorsorformula} to specialise $\phi$ to $\bar{a}_{\widehat{w}_{J}}$ for every $J\in H$.
    Thus, for every $J\in H$ there is for any $J^{\prime}\in H\backslash J$ some $c_{\widehat{w}_{J}}^{ J^{\prime}}\in \bar{a}_{w_{J\cup J^{\prime}}}$ such that for $\Tilde{c}_{\widehat{w}_{J}}=\sum_{J^{\prime}\in H\backslash J}c_{\widehat{w}_{J}}^{ J^{\prime}}$ we have that $\mathfrak{T}_{\Tilde{b}_{\widehat{w}_{J}}+\Tilde{c}_{\widehat{w}_{J}}}$ is realised in $\bar{a}_{\widehat{w}_{J}}$. Finally, we simply choose for any $J^{\prime}\in H\backslash J$ some $k\in J^{\prime}\backslash J$ and consequently define $c_{\widehat{w}_{J}}^{k}$ as the sums of the corresponding $c_{\widehat{w}_{J}}^{ J^{\prime}}$ (and as $0$, if the sum is empty). As $\bar{a}_{w_{J\cup J^{\prime}}}\subseteq \bar{a}_{w_{J\cup\{k\}}}$ we obtain a system as wanted.\\
    \textit{Proof of Claim 2.} Consider the element\[\Tilde{c}=\sum_{J\in H}\Tilde{c}_{\widehat{w}_{J}}=\sum_{J\in H,\,k\not\in J}c_{\widehat{w}_{J}}^{k}.\]
    Then we have that $\Tilde{b}+\Tilde{c}=\sum_{J\in H}\Tilde{b}_{\widehat{w}_{J}}+\Tilde{c}_{\widehat{w}_{J}}$ and since any $\mathfrak{T}_{\Tilde{b}_{\widehat{w}_{J}}+\Tilde{c}_{\widehat{w}_{J}}}$ is realised in $\bar{a}_{\widehat{w}_{j}}$ it follows that
    $\mathfrak{T}_{\Tilde{b}+\Tilde{c}}$ is realised in $\bar{a}_{\mathbf{m}}$.
    Since $\mathfrak{T}_{\Tilde{b}}$ was realised by assumption, it follows that $\mathfrak{T}_{\Tilde{c}}$ is realised in $\bar{a}_{\mathbf{m}}$ as well. This allows us to finally apply the induction hypothesis on $\sum_{J\in H,\,k\not\in J}c_{\widehat{w}_{J}}^{k}$ and thus we obtain a system $g_{\widehat{w}_{J}}^{k}\in \bar{a}_{w_{J,k}}$ as in the claim.
\end{proof}
\begin{notation}
    We write $\wp_{\sigma}(x):=\sigma(x)-x$.
\end{notation}
The following lemma will be preparatory work for the proof of Proposition \ref{torsorclosedyieldsnamalgam}. It can be seen as a generalisation to higher dimensions of the simple fact that if $b_{1}-b_{2}\in F$ and $b_{1},b_{2}$ are independent over some $\sigma$-AS-closed $E$, then $\wp_{\sigma}(b_{1})=\wp_{\sigma}(b_{2})=e\in E$. This already appeared in the proof of Theorem \ref{theorem3amalgamationcharacterisation}.

\begin{lemma}\label{lemmatorsorclosednessdecomposability}
    Let $n\geq 2$. Let $(\bar{a}_{w})_{w\in\mathcal{P}(\mathbf{n})}$ be an independent $n$-system over $E$ such that $E$ is $(n-1)$-$\sigma$-AS-closed for $(\bar{a}_{w})_{w\in\mathcal{P}(\widehat{\mathbf{n}}_{i})}$ for any $1\leq i\leq n$. Let $b\in F\cap \bar{a}_{\mathbf{n}}$ be such that $b=\sum_{1\leq i\leq n}d_{\widehat{w}_{i}}$ for some $d_{\widehat{w}_{i}}\in \bar{a}_{\widehat{w}_{i}}$. Then, we find some system $(e_{\widehat{w}_{i}}^{k})_{1\leq i,k\leq n,\,k\neq i}$ such that $e_{\widehat{w}_{i}}^{k}=-e_{\widehat{w}_{k}}^{i}$, $\mathfrak{T}_{e_{\widehat{w}_{i}}^{k}}$ is realised in $\bar{a}_{\widehat{w}_{i,k}}$ and $\wp_{\sigma}(d_{\widehat{w}_{i}})=\sum_{1\leq k\leq n,\,k\neq i}e_{\widehat{w}_{i}}^{k}$.
\end{lemma}
\begin{proof}
    The proof is by induction with the case $n=2$ already being discussed in the paragraph before the lemma.
    Now, let $n>2$ and $b\in F\cap \bar{a}_{\mathbf{n}}$ be given as in the assumptions. Since $\wp_{\sigma}(b)=0$ we obtain the additive equation $\sum_{1\leq i\leq n}\wp_{\sigma}(d_{\widehat{w}_{i}})=0$. By decomposability (Lemma \ref{lemmadecomposabilityofgeneraladditiveequations}) we have for every $1\leq i\leq n$ that \[\wp_{\sigma}(d_{\widehat{w}_{i}})=\sum_{1\leq k\leq n,\,k\neq i}f_{\widehat{w}_{i}}^{k}\] for $f_{\widehat{w}_{i}}^{k}\in \bar{a}_{\widehat{w}_{i,k}}$. Thus, we can apply $(n-1)$-$\sigma$-AS-closedness to $\wp_{\sigma}(d_{\widehat{w}_{n}})$ and obtain $\wp_{\sigma}(d_{\widehat{w}_{n}})=\sum_{1\leq i\leq n-1}e_{\widehat{w}_{n}}^{i}$ for some $e_{\widehat{w}_{n}}^{i}\in \bar{a}_{\widehat{w}_{i,n}}$ such that any $\mathfrak{T}_{e_{\widehat{w}_{n}}^{i}}$ is realised in $\bar{a}_{\widehat{w}_{i,n}}$, say by $\Tilde{h}_{\widehat{w}_{n}}^{i}$. Consequently, if we set $h_{\widehat{w}_{i}}=d_{\widehat{w}_{i}}+\Tilde{h}_{\widehat{w}_{n}}^{i}$ for any $1\leq i\leq n-1$, we have that $h_{\widehat{w}_{i}}\in \bar{a}_{\widehat{w}_{i}}$ and \[\wp_{\sigma}(h_{\widehat{w}_{i}})=\wp_{\sigma}(d_{\widehat{w}_{i}})+e_{\widehat{w}_{n}}^{i}.\] Combining this we obtain
    \[\sum_{1\leq i\leq n-1}\wp_{\sigma}(h_{\widehat{w}_{i}})=\sum_{1\leq i\leq n-1}\wp_{\sigma}(d_{\widehat{w}_{i}})+e_{\widehat{w}_{n}}^{i}=\sum_{1\leq i\leq n}\wp_{\sigma}(d_{\widehat{w}_{i}})=0.\] 
    Now we apply the induction hypothesis on the equation $0=\sum_{1\leq i\leq n-1}\wp_{\sigma}(h_{\widehat{w}_{i}})$ over $a_{n}$, which is $(n-2)$-$\sigma$-AS-closed for the system $(\Tilde{a}_{w})_{w\in\mathcal{P}(\mathbf{n-1})}$ (where $\Tilde{a}_{w}=a_{w\cup\{n\}}$) by Lemma \ref{lemmantorsorclosednessonestepbelow}. Thus we obtain a system $\left(e_{\widehat{w}_{i}}^{j}\right)_{1\leq i,j\leq n-1,\,i\neq j}$ such that for any $1\leq i,j\leq n-1$ with $i\neq j$ we have
    \[\wp_{\sigma}(d_{\widehat{w}_{i}})+e_{\widehat{w}_{n}}^{i}=\wp_{\sigma}(h_{\widehat{w}_{i}})=\sum_{1\leq j\leq n-1,\,j\neq i}e_{\widehat{w}_{i}}^{j}\;\;\;\text{and}\;\;\;e_{\widehat{w}_{i}}^{j}=-e_{\widehat{w}_{j}}^{i}\]
    and all $\mathfrak{T}_{e_{\widehat{w}_{i}}^{j}}$ are realised in $\bar{a}_{\widehat{w}_{i,j}}$. Thus, if we set for every $1\leq i\leq n-1$ that $e_{\widehat{w}_{i}}^{n}:=-e_{\widehat{w}_{n}}^{i}$, the system $(e_{\widehat{w}_{i}}^{j})_{1\leq i,j\leq n,\,i\neq j}$ fulfills all the requirements.
\end{proof}

\begin{proposition}\label{torsorclosedyieldsnamalgam}
Let $n\geq 3$. The theory $\mathrm{ACFA}^{+}$ has $n$-amalgamation over all $(n-2)$-$\sigma$-AS-closed sets.
\end{proposition}
\begin{proof}
We prove that $(n-2)$-$\sigma$-AS-closedness of $E=\mathrm{acl}_{\sigma}(E)\subseteq\mathcal{M}\models \mathrm{ACFA}^{+}$ implies that any ff-additive equation of height $n$ is ff-decomposable over $E$. This yields the result by Lemma \ref{lemmaffdecomposableimpliesnamalgam}. We proceed by induction. 
For $n=3$, this is given by Theorem \ref{theorem3amalgamationcharacterisation}.\\
Let $n>3$. We fix some independent $n$-system $(\bar{a}_{w})_{w\in\mathcal{P}(\mathbf{n})}$ and some ff-additive equation of height $n$ given by $\sum_{1\leq i\leq n}b_{\widehat{w}_{i}}=0$. Then $\sum_{1\leq i\leq n}b_{\widehat{w}_{i}}=0$ is decomposable (but not yet necessarily ff-decomposable) by Lemma \ref{lemmadecomposabilityofgeneraladditiveequations} with witnessing system $(d_{\widehat{w}_{i}}^{j})_{1\leq i,j\leq n,\,i\neq j}$. In particular for any $1\leq i\leq n$ we have \[(\star)_{1}\;\;\;\;\;\;b_{\widehat{w}_{i}}=\sum_{1\leq j\leq n,\,j\neq i}d_{\widehat{w}_{i}}^{j}.\]
For the moment we work with the equation $b_{\widehat{w}_{n}}=\sum_{1\leq j\leq n-1}d_{\widehat{w}_{n}}^{j}$. It is an equation in the independent $(n-1)$-system $(\bar{a}_{w})_{w\in\mathcal{P}(\mathbf{n-1})}$. Since $E$ is $(n-2)$-$\sigma$-AS-closed, we can apply Lemma \ref{lemmatorsorclosednessdecomposability} and find a system of $\left(e_{\widehat{w}_{i,n}}^{k}\right)$ for $1\leq i,k\leq n-1,\, i\neq k$ and $e_{\widehat{w}_{i,n}}^{k}\in \bar{a}_{\widehat{w}_{i,k,n}}$ such that \[(\star)_{2}\;\;\;\;\;\;\wp_{\sigma}(d_{\widehat{w}_{n}}^{i})=\sum_{1\leq k\leq n-1\,k\neq i}e_{\widehat{w}_{i,n}}^{k}\;\;\;\;\;\;\text{and}\;\;\;\;\;\;(\star)_{3}\;\;\;\;\;\;e_{\widehat{w}_{i,n}}^{k}=-e_{\widehat{w}_{k,n}}^{i}\]hold and moreover for every $e_{\widehat{w}_{i,n}}^{k}$ there is some $f_{\widehat{w}_{i,n}}^{k}\in \mathfrak{T}_{e_{\widehat{w}_{i,n}}^{k}}$ with $f_{\widehat{w}_{i,n}}^{k}\in \bar{a}_{\widehat{w}_{i,k,n}}$. By $(\star)_{3}$ we can assume that the system of the $f_{\widehat{w}_{i,n}}^{k}$ then satisfies for any  $1\leq i,k\leq n-1,\, i\neq k$ that\[(\star)_{4}\;\;\;\;\;\;\wp_{\sigma}(f_{\widehat{w}_{i,n}}^{k})=e_{\widehat{w}_{i,n}}^{k}\;\;\;\;\;\;\text{and}\;\;\;\;\;\;(\star)_{5}\;\;\;\;\;\;f_{\widehat{w}_{i,n}}^{k}=-f_{\widehat{w}_{k,n}}^{i}.\] 
Now we return to the decomposition $(d_{\widehat{w}_{i}}^{j})_{1\leq i,j\leq n,\,j\neq i}$ of our initial equation $\sum_{1\leq i\leq n}b_{\widehat{w}_{i}}=0$ and alter it in the following way:\\
We set for every $1\leq i \leq n-1$: \[c_{\widehat{w}_{n}}^{i}=d_{\widehat{w}_{n}}^{i}-\sum_{1\leq k\leq n-1,\,k\neq i}f_{\widehat{w}_{i,n}}^{k}.\]
Then by $(\star)_{2}$ and $(\star)_{4}$ we have for every $1\leq i\leq n-1$ that
\[\wp_{\sigma}(c_{\widehat{w}_{n}}^{i})=\wp_{\sigma}(d_{\widehat{w}_{n}}^{i})-\sum_{1\leq k\leq n-1,\,k\neq i}\wp_{\sigma}(f_{\widehat{w}_{i,n}}^{k})=\wp_{\sigma}(d_{\widehat{w}_{n}}^{i})-\sum_{1\leq k\leq n-1,\,k\neq i}e_{\widehat{w}_{i,n}}^{k}=0\]
The last line simply says that for any $1\leq i\leq n-1$ we have that\[(\star)_{6}\;\;\;\;\;\;c_{\widehat{w}_{n}}^{i}\in F.\]
Furthermore it follows by $(\star)_{5}$ that
\[\sum_{1\leq i\leq n-1}\sum_{1\leq k\leq n-1,\,k\neq i}f_{\widehat{w}_{i,n}}^{k}=0\]
and consequently we obtain from $(\star)_{1}$ and the definition of the $c_{\widehat{w}_{n}}^{i}$ that
\[(\star)_{7}\;\;\;\;\;\;b_{\widehat{w}_{n}}=\sum_{1\leq i\leq n-1}c_{\widehat{w}_{n}}^{i}.\]
Next, we want to complete $(c_{\widehat{w}_{n}}^{i})_{1\leq i\leq n-1}$ to a system which witnesses that our initial equation $\sum_{1\leq i\leq n}b_{\widehat{w}_{i}}=0$ is ff-decomposable. To do so, we set for any $1\leq i\leq n-1$\[(\star)_{8}\;\;\;\;\;\;g_{\widehat{w}_{i}}:=b_{\widehat{w}_{i}}+c_{\widehat{w}_{n}}^{i}.\] From $(\star)_{7}$ together with the initial equation $\sum_{1\leq i\leq n}b_{\widehat{w}_{i}}=0$ it follows that
\[(\star)_{9}\;\;\;\;\;\;\sum_{1\leq i\leq n-1}g_{\widehat{w}_{i}}=0.\]
But now by $(\star)_{6}$ together with the fact that any $b_{\widehat{w}_{i}}$ was in the fixed field we obtain (from its definition in $(\star)_{8}$) that $g_{\widehat{w}_{i}}\in F$ for any $1\leq i\leq n-1$. Thus, we can consider $(\star)_{9}$ as an ff-additive equation over $a_{n}\supset E$. By Lemma \ref{lemmantorsorclosednessonestepbelow} we can invoke the induction hypothesis on $(\star)_{9}$ and obtain an ff-decomposition $(c_{\widehat{w}_{i}}^{j})_{1\leq i,j\leq n-1,\,j\neq i}$ where for any $1\leq i\leq n-1$ and $j\neq i$ the following holds \[(\star)_{10}\;\;\;\;\;\;g_{\widehat{w}_{i}}=\sum_{1\leq j\leq n-1,\,j\neq i}c_{\widehat{w}_{i}}^{j}\;\;\;\;\;\;\text{and}\;\;\;\;\;\;c_{\widehat{w}_{i}}^{j}=-c_{\widehat{w}_{j}}^{i}.\] We  recall that by definition $b_{\widehat{w}_{i}}=g_{\widehat{w}_{i}}-c_{\widehat{w}_{n}}^{i}$ holds. Hence, if we set $c_{\widehat{w}_{i}}^{n}:=-c_{\widehat{w}_{n}}^{i}$ for every $1\leq i\leq n-1$, then we get for every $1\leq i\leq n$ using $(\star)_{10}$ and the definition of the $g_{\widehat{w}_{i}}$ in $(\star)_{8}$ that \[b_{\widehat{w}_{i}}=\sum_{1\leq j\leq n,\,j\neq i}c_{\widehat{w}_{i}}^{j}.\]
We have now completed the system of fixed field-elements $(c_{\widehat{w}_{i}}^{j})_{1\leq i,j\leq n,\,i\neq j}$ such that it yields ff-decomposabilty of the initial ff-additive equation $\sum_{1\leq i\leq n}b_{\widehat{w}_{i}}=0$. This completes the proof.
\end{proof}

\begin{theorem}\label{theoremnamalgamationovermodels}
In $\mathrm{ACFA}^{+}$ $n$-amalgamation holds over all substructures whose $\mathcal{L}_{\sigma}$-reduct is a model of $\mathrm{ACFA}$.
\end{theorem}
\begin{proof}
    This follows directly from Proposition \ref{torsorclosedyieldsnamalgam} by Proposition \ref{propositionmodelsarentorsorclosed}.
\end{proof}

\section{A counterexample to 4-amalgamation}\label{sectioncounterexample}

The main goal of this section will be to prove that there are sets $E=\mathrm{acl}_{\sigma}(E)\subseteq \mathcal{M}\models \mathrm{ACFA}^{+}$ that are $\sigma$-AS-closed (and thus 3-amalgamation holds over them) but 4-amalgamation does not hold over them. For the author it was surprising to see that being $\sigma$-AS-closed does not suffice to ensure 4-amalgamation, especially, given the fact that in $\mathrm{PF}, \mathrm{ACFA}$ and $\mathrm{PF}^{+}$ the proofs for 3-amalgamation tend to generalise in a relatively straightforward manner to the higher dimensional case. The counter-example for $\mathrm{ACFA}^{+}$ was found after having established the equivalence of 2-$\sigma$-AS-closedness and 4-amalgamation and in this order we will present the result. This equivalence, that is, Lemma \ref{nottwotorsorclosedgivescounterexample}, can be generalised to higher dimensions as discussed in Section \ref{sectionhigherorderequivalence}. However, we give a direct proof because that is much less involved than the general construction.

\begin{lemma}\label{nottwotorsorclosedgivescounterexample}
    Let $E\subseteq\mathcal{M}\models \mathrm{ACFA}^{+}$ be $\sigma$-AS-closed. If $E$ is not $2$-$\sigma$-AS-closed, then 4-amalgamation does not hold over $E$.
\end{lemma}
\begin{proof}
We work in the $\mathcal{L}_{\sigma}$-reduct of $\mathcal{M}$ which is a model of $\mathrm{ACFA}$. We will show that, if $E$ is not 2-$\sigma$-AS-closed, then there is an independent 4-system containing an ff-additive equation that is not ff-decomposable. This suffices by Lemma \ref{lemmanotffdecomposableyieldsnotamalgamation}. The proof has two steps. First we deal with the following degenerate case, and then pass to the general case.\\ 
For the degenerate case assume not only that $E$ is not 2-$\sigma$-AS-closed but that there is an independent 2-system $\{\bar{a}_{12},\bar{a}_{1},\bar{a}_{2}\}$ over $E$ such that for some $b_{1}\in \bar{a}_{1}$ the torsor $\mathfrak{T}_{b_{1}}$ is not realised in $\bar{a}_{1}$ but such that $\mathfrak{T}_{b_{1}}$ is realised in $\bar{a}_{12}$. Note that in particular $\bar{a}_{1}$ is not $\sigma$-AS-closed (and that is witnessed by $\bar{a}_{12}$). The proof now is similar to the one of Fact \ref{theorem3amalgamationcharacterisation} as we will in some sense reconstruct a 3-amalgamation problem over the corner $\bar{a}_{1}$ but such that it lives in a 4-amalgamation problem.\footnote{Note that it is due to our assumption that we can reconstruct a 3-amalgamation problem \textit{inside} a 4-amalgamation problem. Simply taking $\bar{a}_{1}$ not $\sigma$-AS-closed and applying the construction as outlined before Definition \ref{definitionequivrelationEa} does not yield a 4-amalgamation problem because if we take some $\alpha$ with $\sigma(\alpha)-\alpha=b_{1}$, then $\mathrm{acl}_{\sigma}(E\alpha)$ will certainly not be independent from $\bar{a}_{1}$ over $E$.}
Concretely, let $\{\bar{a}_{12},\bar{a}_{1},\bar{a}_{2}\}$ be as above and let $\bar{a}_{13},\bar{a}_{14}$ be copies of $\bar{a}_{12}$ which are chosen such that moreover $\bar{a}_{13}$ is independent from $\bar{a}_{12}$ over $\bar{a}_{1}$ and $\bar{a}_{14}$ is independent from $\bar{a}_{12}\bar{a}_{13}$ over $\bar{a}_{1}$. By taking the corresponding algebraic closures we obtain an independent 4-system $(\bar{a}_{w})_{w\in\mathcal{P}(\mathbf{4})}$ for $\mathrm{ACFA}$.\\
Let $\alpha_{12},\alpha_{13},\alpha_{14}$ be the corresponding realisations of $\mathfrak{T}_{b_{1}}$ such that $\alpha_{1i}\in\bar{a}_{1i}$ for $2\leq i\leq 4$. Let $b_{1ij}:=\alpha_{1i}-\alpha_{1j}$ for $2\leq i,j\leq 4$. We then have $b_{123}+b_{134}+b_{142}=0$. By adding $0$ we obtain an ff-additive equation of height 4 in $(\bar{a}_{w})_{w\in\mathcal{P}(\mathbf{4})}$ which we will show to be not ff-decomposable. To do so, it suffices to show that $b_{123}\neq f_{12}+f_{13}+f_{23}$ for any $f_{ij}\in\bar{a}_{ij}$. Assume the contrary, that is, $b_{123}= f_{12}+f_{13}+f_{23}$ for $f_{ij}\in\bar{a}_{ij}\cap F$. From $b_{123}=\alpha_{12}-\alpha_{13}$ we obtain 
\[(\alpha_{12}-f_{12})-(\alpha_{13}+f_{13})=f_{23}.\]
 Since $\alpha_{12}-f_{12}\in\bar{a}_{12}$ (and similarly for $\alpha_{13}$) it follows by decomposability (Lemma \ref{lemmadecomposabilityofgeneraladditiveequations}) that $f_{23}=r_{2}-r_{3}$ for $r_{2}\in\bar{a}_{2}$ and $r_{3}\in\bar{a}_{3}$. Combining the above equations we obtain
 \[\alpha_{12}-f_{12}-r_{2}=\alpha_{13}+f_{13}-r_{3}\]
We denote the left-hand side by $h_{12}$ (and the right-hand side by $h_{13}$). By independence of $\bar{a}_{12}$ from $\bar{a}_{13}$ over $\bar{a}_{1}$ we get that $h_{12}=h_{13}\in\bar{a}_{1}$. By independence of $\bar{a}_{2}$ from $\bar{a}_{3}$ over $E$ we get $\wp_{\sigma}(r_{2})=\wp_{\sigma}(r_{3})\in E$. Let $e\in E$ with $\wp_{\sigma}(e)=\wp_{\sigma}(r_{2})$ (using that $E$ is $\sigma$-AS-closed). Finally, we obtain (recalling that $f_{12}\in F$) 
\[\wp_{\sigma}(h_{12}+e)=\wp_{\sigma}(\alpha_{12}-f_{12}-r_{2}+e)=\wp_{\sigma}(\alpha_{12})=b_{1}\]
which contradicts the assumption that $\mathfrak{T}_{b_{1}}$ was not realised in $\bar{a}_{1}$.\\
For the rest of the proof, we can thus make the assumption that a 2-system $\{\bar{a}_{12},\bar{a}_{1},\bar{a}_{2}\}$ with the above properties does not exist. We will refer to this assumption as $(\star)$.\\
   Let $\bar{a}_{1},\bar{a}_{2}$ with $b_{1}\in \bar{a}_{1}$ and $b_{2}\in \bar{a}_{2}$ witness that $E$ is not 2-$\sigma$-AS-closed. I.e., $\mathfrak{T}_{b_{1}+b_{2}}$ is realised in $\bar{a}_{12}$ but neither $\mathfrak{T}_{b_{1}}$ nor $\mathfrak{T}_{b_{2}}$ is realised in $\bar{a}_{1}$ or $\bar{a}_{2}$ respectively. 
   As before, by taking independent copies (over $\bar{a}_{1}$) of $\bar{a}_{12}$ we construct an independent 4-system $(\bar{a}_{w})_{w\in\mathcal{P}(\mathbf{4})}$. There we find $b_{3}\in \bar{a}_{3}$ and $b_{4}\in \bar{a}_{4}$ such that $\mathfrak{T}_{b_{1}+b_{3}}$ and $\mathfrak{T}_{b_{1}+b_{4}}$ are realised in $\bar{a}_{13}, \bar{a}_{14}$ respectively but $\mathfrak{T}_{b_{3}}$ and $\mathfrak{T}_{b_{4}}$ are not realised in $\bar{a}_{3},\bar{a}_{4}$ respectively. We will construct an independent $4$-system (in $\mathcal{L}_{\sigma}$) giving rise to an ff-additive equation of height $4$ that is not ff-decomposable. By Lemma \ref{lemmanotffdecomposableyieldsnotamalgamation} this suffices.\\
   Now, since $\mathfrak{T}_{-b_{1}-b_{3}}$ is also realised in $\bar{a}_{13}$ it follows that $\mathfrak{T}_{b_{1}+b_{2}-b_{1}-b_{3}}=\mathfrak{T}_{b_{2}-b_{3}}$ is realised in $\bar{a}_{123}$. But now by $(\star)$ we have that $\mathfrak{T}_{b_{2}-b_{3}}$ is already realised in $\bar{a}_{23}$, because otherwise $\bar{a}_{1}$ and $\bar{a}_{23}$ (together with the element $b_{2}-b_{3}$) would yield a counterexample to $(\star)$ over $E$.\\
   In a similar way, we obtain that $\mathfrak{T}_{b_{2}-b_{4}}$ is realised in $\bar{a}_{24}$ and then consequently that $\mathfrak{T}_{b_{3}-b_{4}}$ is realised in $\bar{a}_{34}$. Now we name the realisations:
   \begin{multicols}{3}
   \begin{itemize}
       \item $c_{12}\models \mathfrak{T}_{b_{1}+b_{2}}$
       \item $c_{13}\models \mathfrak{T}_{b_{1}+b_{3}}$
       \item $c_{14}\models \mathfrak{T}_{b_{1}+b_{4}}$
       \item $c_{23}\models \mathfrak{T}_{b_{2}-b_{3}}$
       \item $c_{24}\models \mathfrak{T}_{b_{2}-b_{4}}$
       \item $c_{34}\models \mathfrak{T}_{b_{3}-b_{4}}$
   \end{itemize}
   \end{multicols}
\noindent Next, we define
\begin{multicols}{2}
\begin{itemize}
    \item $f_{123}=c_{12}-c_{13}-c_{23}$
    \item $f_{124}=-c_{12}+c_{24}+c_{14}$
    \item $f_{134}=c_{13}-c_{34}-c_{14}$
    \item $f_{234}=c_{23}-c_{24}+c_{34}$
\end{itemize}
\end{multicols}
\noindent Then it follows that \[f_{123}+f_{124}+f_{134}+f_{234}=0\]
and, moreover, $f_{ikj}\in F$ for every possible $ikj$. Now it remains to show that $\sum f_{ikj}=0$ is an ff-additive equation that is not ff-decomposable to deduce that $4$-amalgamation does not hold over $E$. Assuming the contrary, we then find a system of elements $e_{ij}\in \bar{a}_{ij}\cap F$ such that
\begin{multicols}{2}
\begin{itemize}
    \item $f_{123}=e_{12}-e_{13}-e_{23}$
    \item $f_{124}=-e_{12}+e_{24}+e_{14}$
    \item $f_{134}=e_{13}-e_{34}-e_{14}$
    \item $f_{123}=e_{23}-e_{24}+e_{34}$ 
\end{itemize}
\end{multicols}
\noindent Now set $d_{ij}=c_{ij}-e_{ij}$. Then, it follows that
\begin{multicols}{2}
\begin{itemize}
    \item $0=d_{12}-d_{13}-d_{23}$
    \item $0=-d_{12}+d_{24}+d_{14}$
    \item $0=d_{13}-d_{34}-d_{14}$
    \item $0=d_{23}-d_{24}+d_{34}$
\end{itemize}
\end{multicols}
\noindent Let us consider the first equation. We obtain by Lemma \ref{lemmadecomposabilityofgeneraladditiveequations} that $d_{12}=g_{1}+g_{2}$ for $g_{1}\in \bar{a}_{1}, g_{2}\in \bar{a}_{2}$. From this it follows that
\[b_{1}+b_{2}=\wp_{\sigma}(c_{12})=\wp_{\sigma}(c_{12}-e_{12})=\wp_{\sigma}(d_{12})=\wp_{\sigma}(g_{1})+\wp_{\sigma}(g_{2})\]
which contradicts the assumption that $b_{1}$ and $b_{2}$ witness that $E$ is not $2$-$\sigma$-AS-closed.
\end{proof}

 The goal of the rest of this section will be to construct a $\sigma$-AS-closed set $E$ that is not $2$-$\sigma$-AS-closed. Consequently, by Lemma \ref{nottwotorsorclosedgivescounterexample}, we will be able to construct $E=\mathrm{acl}_{\sigma}(E)\subseteq \mathcal{M}\models \mathrm{ACFA}^{+}$ such that $3$-amalgamation holds over $E$ but $4$-amalgamation does not. The construction of the counterexample will be \textit{by hand} via concrete computation of solutions of difference equations. It would be very interesting to obtain a more conceptual explanation of the underlying phenomena. Concretely, we will prove the following.
\begin{lemma}\label{theoremcounterexampleexists}
There is some $\sigma$-AS-closed set $E=\mathrm{acl}_{\sigma}(E)\subset M$ and some $e\in E\backslash\{0\}$ such that we can find $a_{1},a_{2}\in M$ independent over $E$ such that
\[\sigma(a_{1})-ea_{1}=e\;\;\;\text{ and }\;\;\;\sigma(a_{2})-\frac{1}{e}a_{2}=\frac{1}{e}\]but such that $\mathfrak{T}_{a_{i}}$ is not realised in $\mathrm{acl}_{\sigma}(Ea_{i})$ for $i=1,2$.
\end{lemma}
Before turning to its proof, we show that Lemma \ref{theoremcounterexampleexists} indeed yields an example as described above.
\begin{theorem}\label{theoremnot4amalagamation}
    There is some $\sigma$-AS-closed set $E=\mathrm{acl}_{\sigma}(E)$ that is not $2$-$\sigma$-AS-closed. Consequently, there is some model of $\mathrm{ACFA}^{+}$ containing an $\sigma$-AS-closed set over which 4-amalgamation does not hold.
\end{theorem}
\begin{proof}
 Let $E,e$ and $a_{1},a_{2}$ be given as by Lemma \ref{theoremcounterexampleexists}. As $\mathfrak{T}_{a_{i}}$ is not realised in $\mathrm{acl}_{\sigma}(Ea_{i})$ for $i=1,2$ it suffices to prove that $\mathfrak{T}_{a_{1}+a_{2}}$ has a solution in $\mathrm{acl}_{\sigma}(Ea_{1}a_{2})$ to obtain that $E$ is not $2$-$\sigma$-AS-closed. We consider the element $a_{1}a_{2}\in \mathrm{acl}_{\sigma}(Ea_{1}a_{2})$ and observe that
 \[\sigma(a_{1}a_{2})-a_{1}a_{2}=(ea_{1}+e)\left(\frac{1}{e}a_{2}+\frac{1}{e}\right)-a_{1}a_{2}=a_{1}+a_{2}+1.\]
 Since $E$ was already $\sigma$-AS-closed, there is some $c\in E$ such that $\sigma(c)-c=1$ and then it follows that $a_{1}a_{2}-c\in \mathfrak{T}_{a_{1}+a_{2}}$ and consequently $E$ is not $2$-$\sigma$-AS-closed.
\end{proof}

Now we collect some preliminary results that will later be useful in the proof of Lemma \ref{theoremcounterexampleexists}.
The first step will be to generalise Lemma \ref{lemmatorsornotrealisedinalgclosure} to $\sigma$-AS equations in the sense of the following definition.

\begin{definition}
    Given $E=\mathrm{acl}_{\sigma}(E)$ we call by a \textit{twisted} $\sigma$-AS equation over $E$ a difference equation in one variable $x$ of the form $\sigma(x)-e_{1}x=e_{2}$ for $e_{1},e_{2}\in E$ where $e_{1}\neq 0$. We call the set it defines a \textit{twisted} $\sigma$-AS set.
\end{definition}
We will now show two lemmas that are both slight generalisations of Lemma \ref{lemmatorsornotrealisedinalgclosure}. The idea of the argument is the same in both cases. However, as the assumptions and calculations differ, the two lemmas will be stated separately. 
\begin{lemma}\label{lemmalineardiffequationnotrealisedinalgclosure}
Let $(K,\sigma)$ be a difference subfield of $M$. Assume that for some $e_{1},e_{2}\in K\backslash\{0\}$ the twisted $\sigma$-AS-equation $\sigma(x)-e_{1}x=e_{2}$ is not realised in $K$, then $\sigma(x)-e_{1}x=e_{2}$ is also not realised in $K^{\mathrm{alg}}$.    
\end{lemma}
\begin{proof}
    Assume the contrary, i.e., $\sigma(b)-e_{1}b=e_{2}$ for some $b\in K^{\mathrm{alg}}\backslash K$. Let $P(X)=X^{n}+c_{1}X^{n-1}+\dots+c_{n}$ be the minimal polynomial of $b$ over $K$. Then the minimal polynomial of $\sigma(b)$ is given by $P_{\sigma}(X)=X^{n}+\sigma(c_{1})X^{n-1}+\dots+\sigma(c_{n})$. Next, since $\sigma(b)=e_{1}b+e_{2}$ we have $P_{\sigma}(e_{1}b+e_{2})=0$ and thus $\frac{1}{e_{1}^{n}}P_{\sigma}(e_{1}b+e_{2})=0$. We write $\frac{1}{e_{1}^{n}}P_{\sigma}(e_{1}X+e_{2})$ as a polynomial in $X$, denoted by $\Tilde{P}(X)=X^{n}+c_{1}^{\prime}X^{n-1}+\dots+c_{n}^{\prime}$. Note that it is indeed monic since we normalised $P_{\sigma}(e_{1}X+e_{2})$ by $\frac{1}{e_{1}^{n}}$.
    Then, we have \[c_{1}^{\prime}=\frac{1}{e_{1}^{n}}\left((ne_{1}^{n-1}e_{2}+e_{1}^{n-1}\sigma(c_{1})\right)=\frac{1}{e_{1}}(ne_{2}+\sigma(c_{1})).\] 
    Since $\Tilde{P}(b)=0$ and $\text{deg}_{X}\Tilde{P}= n$ we have, by uniqueness of the minimal polynomial, $c_{1}^{\prime}=c_{1}$. But from \[c_{1}=\frac{1}{e_{1}}(ne_{2}+\sigma(c_{1}))\] it follows that
\[\sigma\left(-\frac{c_{1}}{n}\right)+e_{1}\frac{c_{1}}{n}=e_{2},\] which contradicts $\sigma(x)-e_{1}x=e_{2}$ not being realised in $K$.
\end{proof}

\begin{lemma}\label{lemmamultiplicativediffequationnotrealisedinalgclosure}
Let $(K,\sigma)$ be a difference subfield of $M$. Assume that for some $e\in K\backslash\{0\}$ the equations $\sigma(x)=e^{z}x$ are not realised in $K\backslash\{0\}$ for any $z\in\mathbb{Z}\backslash\{0\}$, then $\sigma(x)=e^{z}x$ is also not realised in $K^{\mathrm{alg}}\backslash\{0\}$ for any $z\in\mathbb{Z}\backslash\{0\}$.    
\end{lemma}
\begin{proof}
    Assume the contrary, i.e., $\sigma(b)=e^{z}b$ for some $z\in\mathbb{Z}$ and $b\in K^{\mathrm{alg}}\backslash K$. Let $P(X)=X^{n}+c_{1}X^{n-1}+\dots+c_{n}$ be the minimal polynomial of $b$ over $K$. Then the minimal polynomial of $\sigma(b)$ is given by $P_{\sigma}(X)=X^{n}+\sigma(c_{1})X^{n-1}+\dots+\sigma(c_{n})$. Next, since $\sigma(b)=e^{z}b$ we have $P_{\sigma}(e^{z}b)=0$ and thus $\frac{1}{e^{zn}}P_{\sigma}(e^{z}b)=0$. We write $\frac{1}{e^{zn}}P_{\sigma}(e^{z}X)$ as a polynomial in $X$, denoted by $\Tilde{P}(X)=X^{n}+c_{1}^{\prime}X^{n-1}+\dots+c_{n}^{\prime}$ which is then a monic polynomial. \\
    Let $1\leq k\leq n$ be minimal such that $c_{k}\neq 0$. (Note that then $c_{i}^{\prime}=0$ for any $1\leq i<k$ as well.) It follows that 
    \[c_{k}^{\prime}=\frac{1}{e^{zn}}\left(e^{z(n-k)}\sigma(c_{k})\right)=\frac{1}{e^{zk}}\sigma(c_{k}).\] 
    Since $\Tilde{P}(b)=0$ and $\text{deg}_{X}\Tilde{P}= n$ we have by uniqueness of the minimal polynomial that $c_{k}^{\prime}=c_{k}$. But then it follows that \[c_{k}=\frac{1}{e^{zk}}(\sigma(c_{k}))\] 
    which contradicts the assumption that $\sigma(x)-e^{zk}x=0$ is not realised in $K\backslash\{0\}$.
\end{proof}

Let us briefly explain how we will proceed to prove Lemma \ref{theoremcounterexampleexists} and why Lemma \ref{lemmalineardiffequationnotrealisedinalgclosure} and Lemma \ref{lemmamultiplicativediffequationnotrealisedinalgclosure} in this setting are particularly useful. Let $H_{1}(x),H_{2}(x)$ be two twisted $\sigma$-AS equations over $E=\mathrm{acl}_{\sigma}(E)$. We will assume that $a\in M\backslash E$ is a solution of $H_{1}(x)$ and then investigate which condition is necessary for $H_{2}(x)$ having a solution in $\mathrm{acl}_{\sigma}(Ea)$ as well.\\
At this point, applying Lemma \ref{lemmalineardiffequationnotrealisedinalgclosure} and Lemma \ref{lemmamultiplicativediffequationnotrealisedinalgclosure} we can actually reduce to checking whether there is a solution of $H_{2}(x)$ in $\mathrm{cl}_{\sigma}(Ea)$. Now, to deal with this question, we will make use of the assumption that $a$ was a solution of $H_{1}(x)$. Concretely, the main point will be that for any $k\in\mathbb{Z}$ we have that $\sigma^{k}(a)$ is already an element of the (pure) field generated by $Ea$ (since $\sigma(a)=e_{1}a+e_{2}$ for some $e_{1},e_{2}\in E$) and consequently $\mathrm{cl}_{\sigma}(Ea)$ equals the field generated by $Ea$. This allows us to write any element $b\in \mathrm{cl}_{\sigma}(Ea)$ as a fraction of polynomial expressions in $a$ over $E$.\\
Next, we assume that such an element $b$ satisfies $H_{2}(x)$ and we plug the expression of $b$ in terms of $a$ into this equation. Then, by a simple coefficient comparison, we will be able to deduce that certain twisted $\sigma$-AS equations already have to be satisfied in $E$. Using this, we will construct by a simple chain argument a $\sigma$-AS-closed set that avoids a certain set of twisted $\sigma$-AS equations. This will then be enough to prove Lemma \ref{theoremcounterexampleexists}.

\begin{lemma}\label{existencetorsorclosedsetavoidingmultiplicativetorsors}
    There is some $\sigma$-AS-closed set $E\subset M$ and $e\in E\backslash\{0\}$ such that any equation, for any $z\in\mathbb{Z}\backslash\{0\}$ from the following system of equations is \textbf{not} realised in $E\backslash\{0\}$:
    \begin{itemize}
        \item $\sigma(x)/x=e^{z}$
        \item $\sigma(x)-\frac{1}{e}x=\frac{1}{e}$
        \item $\sigma(x)-ex=e$.
    \end{itemize}
\end{lemma}
\begin{proof}
    Let $E=\mathrm{acl}_{\sigma}(E)$ (but not necessarily $\sigma$-AS-closed) with $e\in E\backslash\{0\}$ be given such that the above equations do not have a solution in $E\backslash\{0\}$. Let us briefly explain that we can always find such a set $E$. Take some $D=\mathrm{acl}_{\sigma}(D)\subset M$ and let $g$ be a transformally transcendental element over $D$. Let $\Tilde{D}:=\mathrm{cl}_{\sigma}(Dg)$ and set $E:=\mathrm{acl}(\Tilde{D})$ and $e:=g$, then none of the above equations has a solution in $\Tilde{D}\backslash\{0\}$ and by Lemma \ref{lemmalineardiffequationnotrealisedinalgclosure} and Lemma \ref{lemmamultiplicativediffequationnotrealisedinalgclosure} they do not have a solution in $E\backslash\{0\}$.\\
    We will extend $E$ by a chain construction to a $\sigma$-AS-closed set such that the above equations are still avoided. By compactness it suffices to show that we can realise any $\mathfrak{T}_{f}$ (for $f\in E$) but avoid the above equations in every step. Precisely, we will show that for any $f\in E$ and solution $a\in M\backslash E$ of $\mathfrak{T}_{f}$ none of the above equations holds in $\mathrm{acl}_{\sigma}(Ea)$.\\ 
    First of all, by Lemma \ref{lemmalineardiffequationnotrealisedinalgclosure} and Lemma \ref{lemmamultiplicativediffequationnotrealisedinalgclosure}, it suffices to check that there is no $b\in \mathrm{cl}_{\sigma}(Ea)\backslash E$ being a solution to any of the equations. Since $\sigma(a)=a+f$, for every $b\in \mathrm{cl}_{\sigma}(Ea)\backslash E$ there are $e_{0},\dots,e_{n},\Tilde{e}_{0},\dots,\Tilde{e}_{m}\in E$ with $e_{n},\Tilde{e}_{m}\in E\backslash\{0\}$ such that $b$ can be written as
\[(\star)_{1}\;\;\;b=\frac{\sum_{i=0}^{n}e_{i}a^{i}}{\sum_{j=0}^{m}\Tilde{e}_{j}a^{j}}\;\;\;\text{and consequently}\;\;\;\sigma(b)=\frac{\sum_{i=0}^{n}\sigma(e_{i})(a+f)^{i}}{\sum_{j=0}^{m}\sigma(\Tilde{e}_{j})(a+f)^{j}}.\]
Now assume that $\sigma(b)/b=e^{z}$ for some $z\in\mathbb{Z}\backslash\{0\}$. We plug $(\star)_{1}$ into $\sigma(b)=e^{z}b$ and obtain after multiplying by the denominators
\[\left(\sum_{i=0}^{n}\sigma(e_{i})(a+f)^{i}\right)\left(\sum_{j=0}^{m}\Tilde{e}_{j}a^{j}\right)=e^{z}\left(\sum_{i=0}^{n}e_{i}a^{i}\right)\left(\sum_{j=0}^{m}\sigma(\Tilde{e}_{j})(a+f)^{j}\right)\]
Now we simply compare the coefficients of the highest order terms (in the transcendental element $a$) on both sides and retrieve the equation
\[\sigma(e_{n})\Tilde{e}_{m}=e^{z}e_{n}\sigma(\Tilde{e}_{m})\iff \sigma\left(\frac{e_{n}}{\Tilde{e}_{m}}\right)\left(\frac{e_{n}}{\Tilde{e}_{m}}\right)^{-1}=e^{z}\]
thus contradicting the assumption that $\sigma(x)/x=e^{z}$ does not have a solution in $E\backslash\{0\}$.\\
Next, we assume that $\sigma(b)-eb=e$. Again, we plug $(\star)_{1}$ into the equation $\sigma(b)-eb=e$ and obtain after multiplying by the denominators that

\begin{gather*}
    \left(\sum_{i=0}^{n}\sigma(e_{i})(a+f)^{i}\right)\left(\sum_{j=0}^{m}\Tilde{e}_{j}a^{j}\right)-e\left(\sum_{j=0}^{m}\sigma(\Tilde{e}_{j})(a+f)^{j}\right)\left(\sum_{i=0}^{n}e_{i}a^{i}\right)\\=e\left(\sum_{j=0}^{m}\sigma(\Tilde{e}_{j})(a+f)^{j}\right)\left(\sum_{j=0}^{m}\Tilde{e}_{j}a^{j}\right)
    \end{gather*}
For the equation to hold, since $a$ is not algebraic over $E$, we necessarily have $n\geq m$. We consider the cases of $n>m$ and $n=m$ separately.\\
Assume that $n>m$. In this case, since the highest order term on the right-hand side is of degree $2m$, it follows that we have have cancellation in the highest order terms on the left-hand side, i.e., the following has to hold:
\[\sigma(e_{n})\Tilde{e}_{m}=e\sigma(\Tilde{e}_{m})e_{n}\;\iff\;\sigma\left(\frac{e_{n}}{\Tilde{e}_{m}}\right)\left(\frac{e_{n}}{\Tilde{e}_{m}}\right)^{-1}=e\]
Now this is a contradiction to the assumption that $\sigma(x)/x=e$ does not have a solution in $E\backslash\{0\}$.\\
Next, we assume that $n=m$. In this case the highest order term on both left- and right-hand side are of degree $2m$ and thus the following has to hold:
\[\sigma(e_{n})\Tilde{e}_{m}-e\sigma(\Tilde{e}_{m})e_{n}=e\sigma(\Tilde{e}_{m})\Tilde{e}_{m}\;\iff\; \sigma\left(\frac{e_{n}}{\Tilde{e}_{m}}\right)-e\left(\frac{e_{n}}{\Tilde{e}_{m}}\right)=e\]contradicting the assumption that $\sigma(x)-ex=e$ did not have a solution in $E\backslash\{0\}$.
The proof for the equations $\sigma(x)-\frac{1}{e}x=\frac{1}{e}$ is the same and hence all cases are covered.\end{proof}

Finally, we proceed to prove Lemma \ref{theoremcounterexampleexists}.
\begin{proof}[Proof of Lemma \ref{theoremcounterexampleexists}]
Let $E$ and $e\in E\backslash\{0\}$ be given as by Lemma \ref{existencetorsorclosedsetavoidingmultiplicativetorsors}. Now let $a\in M$ be some realisation of the equation $\sigma(x)-ex=e$. We want to show that there is no $b\in \mathrm{acl}_{\sigma}(Ea)$ such that $b$ realises $\mathfrak{T}_{a}$. The same statement for the equation $\sigma(x)-\frac{1}{e}x=\frac{1}{e}$ then holds as well since replacing $e$ by $\frac{1}{e}$ in the condition of Lemma \ref{existencetorsorclosedsetavoidingmultiplicativetorsors} yields the exact same set of twisted $\sigma$-AS equations that are avoided in $E\backslash\{0\}$.\\
Now, by Lemma \ref{lemmatorsornotrealisedinalgclosure} it suffices to check if there is some $b\in \mathrm{cl}_{\sigma}(Ea)\cap \mathfrak{T}_{a}$. Since $\sigma(a)=e+ea$, it follows that for every $b\in \mathrm{cl}_{\sigma}(Ea)\backslash E$ there are $e_{0},\dots,e_{n},\Tilde{e}_{0},\dots,\Tilde{e}_{m}\in E$ with $e_{n},\Tilde{e}_{m}\in E\backslash\{0\}$ such that $b$ can be written as
\[b=\frac{\sum_{i=0}^{n}e_{i}a^{i}}{\sum_{j=0}^{m}\Tilde{e}_{j}a^{j}}\;\;\;\text{and consequently}\;\;\;\sigma(b)=\frac{\sum_{i=0}^{n}\sigma(e_{i})(ea+e)^{i}}{\sum_{j=0}^{m}\sigma(\Tilde{e}_{j})(ea+e)^{j}}.\]
Now assume that $\sigma(b)-b=a$ holds. Then plugging in the above we obtain after multiplying by the denominators the equality
\begin{gather*}
(\star)\;\;\;\left(\sum_{i=0}^{n}\sigma(e_{i})(ea+e)^{i}\right)\left(\sum_{j=0}^{m}\Tilde{e}_{j}a^{j}\right)-\left(\sum_{i=0}^{n}e_{i}a^{i}\right)\left(\sum_{j=0}^{m}\sigma(\Tilde{e}_{j})(ea+e)^{j}\right)\\
=a\left(\sum_{j=0}^{m}\Tilde{e}_{j}a^{j}\right)\left(\sum_{j=0}^{m}\sigma(\Tilde{e}_{j})(ea+e)^{j}\right)
\end{gather*}
For the equation to hold, using that $a\notin \mathrm{acl}(E)=E$, it is clearly necessary that $n>m$. We will consider the cases $n=m+1$ and $n>m+1$ separately. In both cases we will compare the coefficients of the highest order term to obtain a solution (in $E\backslash\{0\}$) of some twisted $\sigma$-AS equation that was assumed to not have a solution in $E\backslash\{0\}$.\\
Let us start by assuming that $n>m+1$. In this case, since the highest order term on the right-hand side is of degree $2m+1$, it follows that we have cancellation in the highest order terms on the left-hand side, i.e., the following has to hold:
\[\sigma(e_{n})(ea)^{n}\Tilde{e}_{m}a^{m}-e_{n}a^{n}\sigma(\Tilde{e}_{m})(ea)^{m}=0\]
From this we obtain the following \[e^{n}\sigma(e_{n})\Tilde{e}_{m}-e^{m}e_{n}\sigma(\Tilde{e}_{m})=0\;\iff\; e^{n-m}=\frac{e_{n}\sigma(\Tilde{e}_{m})}{\sigma(e_{n})\Tilde{e}_{m}}=\frac{e_{n}}{\Tilde{e}_{m}}\sigma\left(\frac{e_{n}}{\Tilde{e}_{m}}\right)^{-1}.\]
Now this cannot be the case because it would contradict that the equations of the form $\sigma(x)/x=e^{z}$ for $z\in\mathbb{Z}\backslash\{0\}$ do not have a solution in $E\backslash\{0\}$.\\
Now we deal with the case $n=m+1$. In this case the highest order term on both left- and right-hand side are of degree $2m+1$ and thus the following has to hold:
\[e^{n}\sigma(e_{n})\Tilde{e}_{m}-e^{m}e_{n}\sigma(\Tilde{e}_{m})=e^{m}\Tilde{e}_{m}\sigma(\Tilde{e}_{m})\]
We divide this equation by $e^{n}\Tilde{e}_{m}\sigma(\Tilde{e}_{m})$ and obtain
\[\sigma\left(\frac{e_{n}}{\Tilde{e}_{m}}\right)-\frac{1}{e}\frac{e_{n}}{\Tilde{e}_{m}}=\frac{1}{e}.\]
This again can not be true since it would contradict that the equation $\sigma(x)-\frac{1}{e}x=\frac{1}{e}$ does not have a solution in $E\backslash\{0\}$. This completes the proof.
\end{proof}

  


\section{Questions and remarks}\label{sectionquestionsandremarks}

\subsection{Higher amalgamation and $n$-$\sigma$-AS-closedeness}\label{sectionhigherorderequivalence}
As mentioned at the beginning of Section \ref{sectioncounterexample}, the equivalence for any $\sigma$-AS-closed $E$ between being 2-$\sigma$-AS-closed and having 4-amalgamation can be generalised to higher dimensions. One direction is given by Proposition \ref{torsorclosedyieldsnamalgam} and the other is proved in Chapter 4.3 of the author's PhD thesis \cite{ludwig:tel-05236078}. We decided to not give the whole argument here due to its technicality and lack of application. The main idea is, given some $\mathcal{L}_{\sigma}(E)$-type $p_{\mathbf{n}}(\bar{x}_{\mathbf{n}})$ that witnesses that $E$ is not $n$-$\sigma$-AS-closed, to construct an ff-additive equation of height $(n+2)$ that is not ff-decomposable. To do so one constructs a type such that several prescribed sub-tuples of any of its realisations are models of $p_{\mathbf{n}}$ and then using a combinatorial argument one reorders the sub-tuples in order to obtain an independent $(n+2)$ system that witnesses the existence of the ff-additive equation that is not decomposable.

\subsection{Definable groupoids}
In \cite{grupoidshrushovski} a connection between definable groupoids, internal covers and (the failure of) $4$-amalgamation was established which was then further studied in a series of papers including \cite{grupoidscoversgoodrickkolesnikov},\cite{Goodrick2010AmalgamationFA} and \cite{homologygroupsoftypesgoodrickkimkolesnikov}. Whereas the notion of a definable groupoid can be defined in any theory, the main result of this type in \cite{grupoidshrushovski} states that in a stable theory 4-amalgamation holds if and only if all definable finitary groupoids are eliminable. In \cite{grupoidscoversgoodrickkolesnikov} the slightly stronger notion of retractability of a type-definable groupoid was introduced and linked to the above. When turning from stable to simple theories the above correspondence does not go through in general as shown in \cite{Goodrick2010AmalgamationFA}. However, it was conjectured in \cite{Goodrick2010AmalgamationFA} that failure of 4-amalgamation comes from some (in a suitable sense) interpretable non-eliminable groupoid or from an interpretable tetrahedron-free hypergraph.\\
In $\mathrm{ACFA}^{+}$ as we work in continuous logic it seems necessary to work with type-definable groupoids (and compact instead of finitary) and it would be interesting to see whether the above described analogy extends. A concrete question, that takes into account the equivalence of $4$-amalgamation and $2$-$\sigma$-AS-closedness, would be: \textit{Can we associate to every $E$ which is $\sigma$-AS-closed but not $2$-$\sigma$-AS-closed some (type)-definable groupoid that is not eliminable? How does such a groupoid look like?}

\subsection{Higher amalgamation in $\mathrm{ACFA}^{+,\times}$} In \cite[7.2]{ludwig2025modeltheorydifferencefields} it was already proposed to study the theory $\mathrm{ACFA}^{\times}$ and $\mathrm{ACFA}^{+,\times}$ consisting of the expansions of $\mathrm{ACFA}$ (resp. $\mathrm{ACFA}^{+}$) by a sufficiently generic multiplicative character (added as a continuous logic predicate) on the fixed field using the corresponding expansions of $\mathrm{PF}$ (resp. $\mathrm{PF}^{+}$) introduced in \cite{ludwig2025pseudofinitefieldsadditivemultiplicative}. While the results of \cite{ludwig2025modeltheorydifferencefields} are likely transferable to the multiplicative context (see \cite[7.2]{ludwig2025modeltheorydifferencefields}), for higher amalgamation this is not necessarily the case and a more thorough analysis would be required. Both, the counterexample to 4-amalgamation which uses the action of the multiplicative on the additive group, as well as Lemma \ref{lemmatorsornotrealisedinalgclosure}, the crucial ingredient to the proof of higher amalgamation over models, might not have a straightforward multiplicative analogue. In any way it appears an intriguing question to study higher amalgamation for the purely multiplicative case as well as for $\mathrm{ACFA}^{+,\times}$.

\bibliography{bibpfplustimes}
\bibliographystyle{abbrv}
\end{document}